# On the finite element analysis of functionally graded sandwich curved beams via a new refined higher shear deformation theory

Mohamed-Ouejdi Belarbi[*1], Mohammed Sid Ahmed Houari[2], Hicham Hirane[1], Ahmed Amine Daikh[2], Stéphane Pierre Alain Bordas[3,4]

[1]*Laboratoire de Génie Energétique et Matériaux, LGEM, Université de Biskra, B.P. 145, R.P. 07000, Biskra, Algeria*
[2] *Laboratoire d'Etude des Structures et de Mécanique des Matériaux, Département de Génie Civil, Faculté des Sciences et de la Technologie, Université Mustapha Stambouli B.P. 305, R.P. 29000 Mascara, Algérie.*
[3]*Institute of Research and Development, Duy Tan University, K7/25 Quang Trung, Danang, Vietnam. Department of Engineering*
[4]*Institute of Computational Engineering, University of Luxembourg, Maison du Nombre, 6 Avenue de la Fonte, 4364 Esch-sur-Alzette, Luxembourg.*

**Abstract**

In the present paper, a new parabolic shear deformation beam theory is developed and applied to investigate the bending behavior of functionally graded (FG) sandwich curved beam. The present theory is exploited to satisfy parabolic variation of shear stress distribution along the thickness direction thereby obviating the use of any shear correction factors. The material properties of FG sandwich beam change continuously from one surface to another according to a power-law function. Three common configurations of FG beams are used for the study, namely: (a) single layer FG beam; (b) sandwich beam with FG face sheets and homogeneous core and (c) sandwich beams with homogeneous face sheets and FG core. The governing equations derived herein are solved by employing the finite element method using a two-noded beam element, developed for this purpose. The robustness and reliability of the developed finite element model are demonstrated by comparing its results with those available by other researchers in existing literature. The comparison studies show that the proposed model is: (a) accurate and comparable with the literature; b) of fast rate of convergence to the reference solution; c) excellent in terms of numerical stability and d) valid for FG sandwich curved beams. Moreover, comprehensive numerical results are presented and discussed in detail to investigate the effects of volume fraction index, radius of curvature, material distributions, length-to-thickness ratio, face-to-core- thickness ratio, loadings and boundary conditions on the static response of FG curved sandwich beam. New referential results are reported which will be serve as a benchmark for future research.

[*] Corresponding author, E-mail: mo.belarbi@univ-biskra.dz

# 1. Introduction

Functionally graded materials (FGM) are a relatively new technology which are increasingly being used in components exposed to high temperature gradients [1]. Due to the advent of More Electric Aircrafts (MEA), the emphasis on novel materials like FGMs are even more due to the need of advanced Thermal Management Systems (TMS) [2-3]. Due to the gradual variation of material properties through the thickness, FGMs are used to solve a plethora of engineering problems in marine, automotive, mechanical and civil engineering sectors [4-5]. As a result, a number of studies have been performed to analyze the static, vibration, and buckling response of FG sandwich structures (Plates, Beams, Shells….) [6-20]. A curved composite beam provides additional flexibility due to its geometry and has found its way in many engineering applications (bridges, aircraft, spacecraft, etc…). Various beam theories have been proposed to predict the response of such beams. They can be classified into three groups, namely: classical Euler–Bernoulli beam theory (CBT), first-order shear deformation beam theory (FSDT) and higher-order beam theories (HSDTs). Since the CBT does not incorporate transverse shear deformation effects, its application is limited to very long beams only. The FSDT surmounts this problem by taking into account this effect and gives acceptable results for moderately short and long beams. However, the FSDT needs a shear correction factor which is complex to determine due to its dependence on the geometry, material properties and boundary conditions of each specific problem [21]. The HSDTs on the other hand do not require any shear correction factor and have been found to compute deflections and stresses more accurately. These theories include higher-order terms in the approximation of the in-plane displacement fields and satisfy zero shear stress conditions at the top and the bottom surfaces of the beams.

Many computational models, both analytical and numerical, have been proposed and developed over the years. Frostig et al. [22, 23] developed a new higher order theory based on variational principles to analyze the bending behavior of sandwich beam with transversely flexible core. Venkataraman and Sankar [24] investigated the bending behavior of sandwich beam having a functionally graded (FG) core using an elasticity solution. In their study, it was assumed that the Young's modulus of the core is varied exponentially through the thickness. Sankar [25] presented a three-dimensional (3D) elasticity solution for simply supported FG beams subjected to transverse loads. Daouadji et al. [26] used a partial differential equation to study the static problem of a cantilever FG beam subjected to linearly distributed load. Apetre et al. [27] carried out a comparaison study to investigate the static behavior of a sandwich beam with FG core. In a later study, Şimşek [28] analyzed the static behavior of a simply supported FG beam subjected to a uniformly distributed load using the Ritz method. Li et al. [29] investigated the static bending and dynamic response of FG beams using the HSDT. Based on various HSDTs, Thai and Vo [30] employed Navier technique to obtain the analytical solution of a simply supported FG beam. Larbi et al. [31] proposed an efficient hyperbolic shear deformation beam theory based on the neutral surface position for bending and free vibration analysis of simply supported FG beams. Nguyen and

Nguyen [32] presented a new inverse trigonometric shear deformation theory for static, free vibration and buckling responses of FG sandwich beams. Karamanlı [33] investigated the static behavior of two-directional FG sandwich beam under various boundary conditions by using a quasi-3D shear deformation theory. More recently, a unified five unknown shear deformation theory is developed by Sayyad and Ghugal [34] to analyze the bending response of FG sandwich beams and plates with softcore and hardcore.

From the previous literature review, the majority of researchers used analytical models to investigate the behavior of FG sandwich beam. However, the analytical solutions (2D/3D) were limited to simple geometries, certain types of gradation of material properties (e.g., exponential or power law distribution), special loading cases and specific types of boundary conditions [35, 36]. Therefore, the numerical methods can serve as a better choice analyze the complex behavior of FGM structures. Among them, the finite element method (FEM) is the most popular one. The FEM has several advantages in terms of ease in implementation of complex loading, arbitrary grading properties, varying boundary conditions and ease of solutions process [37-40]. Chakraborty et al. [41] developed a new beam finite element based on the FSDT to investigate the static, wave propagation and free vibration behavior of FG beam. Kadoli et al. [42] presented a new finite element beam formulation based on the TSDT to analyze the static behavior of FG beam under ambient temperature. Based on the same theory, Vo et al. [43] developed a two-noded Hermite-cubic beam element with ten degree of freedom (DOFs) to study the static and vibration responses of FG beam. Later, beam finite element based on 1D Carrera Unified Formulation (CUF) is presented by Filippi et al. [44] to examine the static problem of FG beam using various displacement theories (trigonometric, polynomial, exponential and miscellaneous expansions). Vo et al. [45, 46] developed a two-noded $C^1$ beam element with six DOFs per node and used Navier solutions to determine the displacement, stresses, natural frequencies and critical buckling loads of FG sandwich beam by using a quasi-3D polynomial shear deformation theory. With the aid of zig-zag theory (ZIGT), Khan et al. [47] constructed a two noded beam element having four DOFs at each node for the static and free vibration analysis of FG beam. The authors used linear Lagrange interpolation function for the axial displacement and cubic hermite interpolation for the deflection. Similarly, a new Hermite-Lagrangian finite element with saven DOFs per node is developed by Yarasca et al. [48] to study the bending analysis of FG sandwich beams. The element is formulated based on the hybrid quasi-3D shear deformation theory. Frikha et al. [49] developed a new two noded mixed finite element with eight DOFs for FG beams based on the HSDT. Jing et al. [50] studied the the static bending and free vibration characteristics of FG beam by using a new formulation based on combination of cell-center finite volume method and Timoshenko beam theory. Recently, based on the Quasi-3D HSDT, Nguyen et al. [51] developed an efficient two noded beam element having five DOFs per node to study the static bending of FG beams under various boundary condition. Koutoati et al. [52] sudied the static and free vibration behavior of FGM sandwich beam using three finite element models based on CPT, FSDT and

HSDT. It is concluded that the axial bending and shear coupling affect the response of the FGM sandwich beam in both statics and dynamics. Li et al. [53] established a new mixed finite element beam formulation based on the HSDT for accurate analysis of the FG sandwich beam with introducing the stress equilibrium condition. More recently, the same authors [54] formulated a new three node beam element with ten DOFs using HSDT to determine the displacement, stresses, and critical buckling loads of FG beams with arbitrary material distribution through thickness. In the same year, Katili et al. [55 developed a two noded Hermitian finite element having four DOFs per node to solve the static and free vibration problems of FG beam. The formulation of this new element is based on the unified and integrated approach of Timoshenko beam theory.

Based on the aforementioned review, it appears that literature on the analysis of static response of FG single layer and straight sandwich beams is plenty. Unfortunately, there is limited work available in the literature for bending analysis of curved beams. Fereidoon et al. [56] studied the bending behavior of a simply supported curved sandwich beam with FG core. They used the CBT to model the thin face-sheets and the HSDT for the core layer. Kurtaran [57] used the generalized differential quadrature method (GDQM) to study the bending and transient behaviors of moderately thick FG deep curved beam. Based on the variational iterational method, Eroglu [58] investigated the large deflection of FG curved planar beams. The flexural response of curved multilayered beam with constant curvature was studied by Thurnherr et al. [59] using higher order beam model. More recently, Sayyad and Ghugal [60] used a quasi-3D sinusoidal shear deformation theory (SSDT) to investigate the static behavior of simply suported symmetric curved sandwich beam with FG face sheets. Similarly, Avhad and Sayyad [61] studied the static deformation of simply suported FG curved sandwich beams using a new polynomial fifth order shear deformation theory. Lezgy-Nazargah [62] develoepd a new finite element model with thirteen DOFs for the static analysis of curved thin-walled beams. The element is formulated by using a global–local layered beam theory.

As far as the authors are aware, there is currently no publication available that explains the bending behavior of symmetric and non-symmetric curved FG sandwich beam with various boundary conditions and arbitrary FG material distribution using finite element method. Therefore, the development of an efficient beam element is necessary to analyze the complex behavior of FG sandwich beam. In the present work, an efficient finite element model to investigate the bending behavior of FG curved sandwich beam has been developed. This new element is formulated based on the recently proposed parabolic shear deformation theory. The present theory provides a parabolic distribution of transverse shear stress across the thickness and satisfies the zero traction boundary conditions on the top and the bottom surfaces obviating the use of any shear correction factors. The governing equations are derived using the virtual work principle and the material properties are varied according to a power-law function. The efficiency of the proposed beam element is demonstrated for symmetric and non-symmetric FG curved sandwich beams with arbitrary FG material distribution, various boundary conditions, face-to-core-thickness ratio, length-to-thickness ratio and volume fraction index. The results are compared with those obtained using the refined analytical solutions and other finite element models availaible in the literature. Finally, several additional results are obtained which will potentially serve as a benchmark for the future investigation.

## 2. Theoretical formulation

### 2.1 Geometrical configuration

A rectangular curved FG sandwich beam with uniform thickness $h$, length $L$ and width $b$ is considered (Fig. 1). The mid-plane of the beam ($z=0$) is considered as the reference plane. The top and bottom surfaces of the plate are at $z=\pm h/2$, and the edges of the plate are parallel to the $x$-axis. Three types of FG beams are studied: a single-layered FG beam (type A); a sandwich beam with FGM face sheets with a homogeneous core (type B); sandwich beam with homogeneous face sheets and FGM core (type C).

### 2.2 Material proprieties

The effective material properties of the FG beam are assumed to vary smoothly across the thickness direction according to a power-law function. They are calculated by using the following rule of mixture [63, 64]:

$$P(z) = P_m + (P_c - P_m)V(z) \tag{1}$$

Note that $P_m$ and $P_c$ are, respectively, the corresponding properties of the metal and the ceramic, $V^{(n)}$ is the volume fraction of each layer $n$ ($n$ = 1, 2, 3). For simplicity, Poisson's ratio of the FG beam is assumed to be constant through thickness in this analysis.

#### 2.2.1 Type A: Isotropic FG beams

The beam is graded from a mixture of metal and ceramic, in which the composition is varied from the bottom surface to the top surface (Fig. 1b). The volume fraction of the FG beam varies along the thickness direction via a power-law function as follows:

$$V_c(z) = \left(\frac{2z+h}{2h}\right)^p, \quad z \in [-h/2, h/2] \tag{2}$$

The parameter "$p$" is the volume fraction index ($0 \leq p \leq +\infty$) that allows the user to define gradation of material properties through the thickness direction. The value of "$p$" equal to 0 and $+\infty$ represents a fully ceramic and fully metal beam, respectively.

#### 2.2.2 Functionally graded sandwich beams

The sandwich beam is composed of three layers. The vertical coordinates of the bottom surface, the two interfaces, and the top surface are denoted by $h_1 = -h/2$, $h_2$, $h_3$, $h_4 = h/2$, respectively. For the brevity, the ratio of the thickness of each layer from bottom to top is denoted by the combination of three numbers, i.e. ''1-0-1'', ''2-1-2'', ''3-4-3'' and so on. As shown in Figs. 1c, d, two types of FG sandwich beam are considered:

- Type B: FG face sheets and homogeneous core.
- Type C: Homogeneous face sheets and FG core.

**2.2.2.1 Type B: Sandwich beams with FG face sheets and homogeneous core**

As shown in Fig. 1c, the top and bottom face sheets are graded from metal to ceramic whereas the core is made of fully ceramic. The volume fraction of the FGMs is assumed to obey a power-law function along the thickness direction:

$$V^{(1)} = \left(\frac{z-h_1}{h_2-h_1}\right)^p, \quad z \in [h_1, h_2] \tag{3a}$$

$$V^{(2)} = 1, \quad z \in [h_2, h_3] \tag{3b}$$

$$V^{(3)} = \left(\frac{z-h_4}{h_3-h_4}\right)^p, \quad z \in [h_3, h_4] \tag{3c}$$

**2.2.3 Type C: sandwich beams with homogeneous face sheets and FG core.**

The volume fraction of the FG core is assumed to obey a power-law function along the thickness direction:

$$V^{(1)} = 0, \quad z \in [h_1, h_2] \tag{4a}$$

$$V^{(2)} = \left(\frac{z-h_2}{h_3-h_2}\right)^k, \quad z \in [h_2, h_3] \tag{4b}$$

$$V^{(3)} = 1, \quad z \in [h_3, h_4] \tag{4c}$$

Fig.2 shows the through thickness variation of the volume fraction function of the mentioned three cases of FG beam for various values of the power law index $p$.

**2.3 Kinematics of the present theory:**

**2.3.1 Displacement field**

A novel quasi-2D parabolic shear deformation beam theory for FG curved beam considering the transverse shear deformation is adopted in this study. The displacement field of the proposed theory is chosen based on the following assumptions:

**(1)** The axial displacement consists of extension, bending and shear components;

**(2)** The bending component of axial displacement is similar to that given by the Euler–Bernoulli beam theory;

**(3)** The shear component of axial displacement gives rise to the parabolic variation of shear strain and hence to shear stress through the thickness of the beam in such a way that shear stress vanishes on the top and bottom surfaces.

Based on the assumptions made above, the displacement field can be obtained as:

$$u(x,z) = \left(1 + \frac{z}{R}\right)u_0(x) - z\frac{\partial w_0}{\partial x} + f(z)\phi_x(x)$$

$$w(x,z) = w_0(x)$$

(5)

where $R$ is the radius of curvature; $u_0(x)$ is the curvilinear axial displacement; $w_0(x)$ is the transverse displacement of the mid-line points of the beam; $\phi_x$ is the rotation of the cross section of the beam at the neutral axis due to transverse shear deformation.

A new parabolic shear deformation beam function is used [65]:

$$f(z) = z\left(1 - \frac{3}{2}\left(\frac{z}{h}\right)^2 + \frac{2}{5}\left(\frac{z}{h}\right)^4\right)$$

(6)

The function $g(z)$ is given as follows:

$$g(z) = f'(z)$$

(7)

### 2.3.2 The strain field

The strain components deduced classically with respect to the curvilinear covariant basis vector are given here as:

$$\varepsilon_x = \frac{\partial u_0}{\partial x} - z\frac{\partial^2 w_0}{\partial^2 x} + f(z)\frac{\partial \phi_x}{\partial x} + \frac{w_0}{R}$$

$$\gamma_{xz} = g(z)\phi_x$$

(8)

Rewrite the strain components in the short form as follows:

$$\varepsilon_x = \varepsilon_x^0 + z\varepsilon_x^1 + f(z)\varepsilon_x^2$$

$$\gamma_{xz} = g(z)\gamma_{xz}^0$$

(9)

where

$$\varepsilon_x^0 = \left(\frac{\partial u_0}{\partial x} + \frac{w_0}{R}\right), \quad \varepsilon_x^1 = \frac{\partial^2 w_0}{\partial^2 x}, \quad \varepsilon_x^2 = \frac{\partial \phi_x}{\partial x}, \quad \gamma_{xz}^0 = \phi_x$$

(10)

The curved beam is made of FG materials. Thus, the constitutive relations between the stress and the strain are given as follows:

$$\begin{Bmatrix} \sigma_x \\ \tau_{xz} \end{Bmatrix}^{(n)} = \begin{bmatrix} C_{11} & 0 \\ 0 & C_{55} \end{bmatrix}^{(n)} \begin{Bmatrix} \varepsilon_x \\ \gamma_{xz} \end{Bmatrix}^{(n)}$$

(11)

where $\sigma_x$ and $\tau_{xz}$ are the axial and transverse shear stresses. $C_{ij}$ are the stiffness coefficients correlated with the engineering constants as follows:

$$C_{11} = E(z)^n \text{ and } C_{55} = \frac{E(z)^n}{2(1+\nu)} \tag{12}$$

## 3. Principle of Virtual Work

The total virtual work principle, considering the static analysis, can be given as:

$$0 = \int_{-h/2}^{h/2} \int_x \delta U - \delta V \tag{13}$$

$$0 = \int_{-h/2}^{h/2} \int_x \delta \varepsilon^T \sigma \, dxdz - \int_x q \delta w dx \tag{14}$$

where $\delta U$ is the internal virtual work and $\delta V$ denote the external virtual work.

Using the strain expression in Eq. (8), the internal virtual work performed by the axial and tangential stresses can be derived as below:

$$\begin{aligned}
\delta U &= \int_0^L \int_{-\frac{h}{2}}^{\frac{h}{2}} \left( \sigma_x \delta \varepsilon_x + \tau_{xz} \delta \gamma_{xz} \right) dzdx \\
&= \int_0^L \left( N_x \frac{d\delta u_0}{dx} - M_x \frac{d^2 \delta w_0}{dx^2} + S_x \frac{d\delta \phi_x}{dx} + N_x \frac{\delta w_0}{R} + Q_{xz} \left[ \delta \phi_x \right] \right) dx
\end{aligned} \tag{15}$$

where $N_x$, $M_x$, $S_x$ and $Q_{xz}$ are, respectively, the axial force, bending moment, shear moment, and shear force. They are defined by:

$$(N_x, M_x, S_x) = \sum_{n=1}^{3} \int_{h_n}^{h_{n+1}} (1, z, f(z)) \sigma_x dz \tag{16a}$$

$$Q_{xz} = \sum_{n=1}^{3} \int_{h_n}^{h_{n+1}} g(z) \tau_{xz} dz \tag{16b}$$

By substituting Eqs. (11) and (8) into Eq. (16), the final expressions for the stress resultants are given as:

$$\begin{aligned}
N_x &= A_{11} \left( \frac{du_0}{dx} + \frac{w_0}{R} \right) - B_{11} \frac{d^2 w_0}{dx^2} + B_{11}^s \frac{d\phi_x}{dx} \\
M_x &= B_{11} \left( \frac{du_0}{dx} + \frac{w_0}{R} \right) - D_{11} \frac{d^2 w_0}{dx^2} + D_{11}^s \frac{d\phi_x}{dx} \\
S_x &= B_{11}^s \left( \frac{du_0}{dx} + \frac{w_0}{R} \right) - D_{11}^s \frac{d^2 w_0}{dx^2} + H_{11}^s \frac{d\phi_x}{dx} \\
Q_{xz} &= A_{55}^s \phi_x
\end{aligned} \tag{17}$$

Rewrite Equation (17) into the matrix form as below:

$$\begin{Bmatrix} N_x \\ M_x \\ S_x \\ Q_{xz} \end{Bmatrix} = \begin{bmatrix} A_{11} & B_{11} & B_{11}^s & 0 \\ B_{11} & D_{11} & D_{11}^s & 0 \\ B_{11}^s & D_{11}^s & H_{11}^s & 0 \\ 0 & 0 & 0 & A_{55}^s \end{bmatrix} \begin{Bmatrix} \varepsilon_x^0 \\ \varepsilon_x^1 \\ \varepsilon_x^2 \\ \gamma_{xz}^0 \end{Bmatrix} \quad (18)$$

where the cross-sectional rigidities are expressed as:

$$\left(A_{11}, B_{11}, D_{11}, B_{11}^s, D_{11}^s, H_{11}^s\right) = \sum_{n=1}^{3} \int_{h_n}^{h_{n+1}} C_{11}\left(1, z, z^2, f(z), z\, f(z), f^2(z)\right) dz. \quad (19a)$$

$$A_{55}^s = \sum_{n=1}^{3} \int_{h_n}^{h_{n+1}} C_{33} g(z)^2 dz. \quad (19b)$$

The external virtual work carried out by the distributed load $q(x)$ can be given as:

$$\delta V = \int_0^L q(x) \delta w \, dx \quad (20)$$

The following weak statement is obtained by using the virtual work principle:

$$0 = \int_0^L \left( N_x \frac{d\delta u_0}{dx} - M_x \frac{d^2 \delta w_0}{dx^2} + S_x \frac{d\delta \phi_x}{dx} + N_x \frac{\delta w_0}{R} + Q_{xz}\left[\delta \phi_x\right] - q(x)\delta w \right) dx \quad (21)$$

By substituting Eqs. (18) into Eq. (21), the final expression of virtual work principle can be written as follows:

$$0 = \int_0^L (\delta\varepsilon^{0^T} A_{11} \varepsilon^0 + \delta\varepsilon^{0^T} B_{11} \varepsilon^1 + \delta\varepsilon^{0^T} B_{11}^s \varepsilon^2 + \delta\varepsilon^{1^T} B_{11} \varepsilon^0 + \delta\varepsilon^{1^T} D_{11} \varepsilon^1 + \delta\varepsilon^{1^T} D_{11}^s \varepsilon^2 +$$

$$\delta\varepsilon^{2^T} B_{11}^s \varepsilon^0 + \delta\varepsilon^{2^T} D_{11}^s \varepsilon^1 + \delta\varepsilon^{2^T} H_{11}^s \varepsilon^2 + \delta\gamma^{0^T} A_{55}^s \gamma^0) dx - q \int_0^L \delta w \, dx \quad (22)$$

For the static analysis, the governing equations associated with the present parabolic shear deformation beam theory are obtained by integrating Eq. (13) by parts. Thus, the following equilibrium equations are obtained by collecting the coefficients of $\delta u_0, \delta w_0, \delta \phi_x$ and equating them with zero.

$$\begin{aligned} \delta u_0: & \quad \frac{dN_x}{dx} = 0 \\ \delta w_0: & \quad \frac{d^2 M_x}{dx^2} - \frac{N_x}{R} + q(x) = 0 \\ \delta \phi_x: & \quad \frac{dS_x}{dx} - Q_{xz} = 0 \end{aligned} \quad (23)$$

## 4. Finite element formulation

In the present study, a two noded beam element having four degrees of freedom (DOFs) at each node is originally developed (Fig. 3) to analyse the static behavior of FG curved sandwich beam.

This new element is formulated based on recently proposed refined higher shear deformation beam theory. The unknown $u_0$ and $\varphi_x$ have been interpolated using $C^0$ linear Lagrange interpolation function ($N$) while the unknown $w_0$ has been interpolated using $C^1$ Hermite cubic interpolation function ($\overline{N}$). The displacement vector corresponding to node $i$ ($i = 1$ to $2$) is given as:

$$d = \left\{ u_{0i} \quad w_{0i} \quad \frac{\partial w_{0i}}{\partial x} \quad \varphi_{xi} \right\}^T \tag{24}$$

The generalized displacements within an element are given as

$$u_0 = N\, u_0^e, \quad \varphi_x = N\, \varphi_x^e, \quad w_0 = \overline{N}\, w_0^e \tag{25}$$

where the nodal displacements of an element $u_0^e$, $\varphi_x^e$ and $w_0^e$ are expressed as

$$u_0^e = \begin{bmatrix} u_{0_1} & u_{0_2} \end{bmatrix}^T, \quad \varphi_x^e = \begin{bmatrix} \varphi_{x_1} & \varphi_{x_2} \end{bmatrix}^T$$
$$w_0^e = \begin{bmatrix} w_{0_1} & w_{0,x_1} & w_{0_2} & w_{0,x_2} \end{bmatrix}^T \tag{26}$$

The classical interpolation functions are defined as:

$$N = \begin{bmatrix} N_1 & N_2 \end{bmatrix}, \quad \overline{N} = \begin{bmatrix} \overline{N}_1 & \overline{N}_2 & \overline{N}_3 & \overline{N}_4 \end{bmatrix} \tag{27}$$

where

$$N_1 = 1 - \frac{x}{L}, \qquad N_2 = \frac{x}{L}$$
$$\overline{N}_1 = 1 - \frac{3x^2}{L^2} + \frac{2x^3}{L^3}, \qquad \overline{N}_2 = x - \frac{2x^2}{L} + \frac{x^3}{L^2} \tag{28}$$
$$\overline{N}_3 = \frac{3x^2}{L^2} - \frac{2x^3}{L^3}, \qquad \overline{N}_4 = -\frac{x^2}{L} + \frac{x^3}{L^2}$$

Substituting Eq. (27) into the generalized strain vectors in Eq. (8) gives:

$$\{\varepsilon^0\}_{(e)} = [B_0]_{(e)} \{\delta_i\}_{(e)}, \quad \{\varepsilon^1\}_{(e)} = [B_1]_{(e)} \{\delta_i\}_{(e)}$$
$$\{\varepsilon^2\}_{(e)} = [B_2]_{(e)} \{\delta_i\}_{(e)}, \quad \{\gamma^0\}_{(e)} = [B_s]_{(e)} \{\delta_i\}_{(e)} \tag{29}$$

where the components of strain–displacement matrices $[B_i]$ for FG sandwich curved beam element are given by:

$$[B_0] = \begin{bmatrix} \dfrac{dN_1}{dx} & \dfrac{1}{R}\overline{N}_1 & \dfrac{1}{R}\overline{N}_2 & 0 & \dfrac{dN_2}{dx} & \dfrac{1}{R}\overline{N}_3 & \dfrac{1}{R}\overline{N}_4 & 0 \\ 0 & 0 & 0 & 0 & 0 & 0 & 0 & 0 \end{bmatrix}$$

$$[B_1] = \begin{bmatrix} 0 & -\dfrac{d^2\overline{N}_1}{dx^2} & -\dfrac{d^2\overline{N}_2}{dx^2} & 0 & 0 & -\dfrac{d^2\overline{N}_3}{dx^2} & -\dfrac{d^2\overline{N}_4}{dx^2} & 0 \\ 0 & 0 & 0 & 0 & 0 & 0 & 0 & 0 \end{bmatrix} \quad (30)$$

$$[B_2] = \begin{bmatrix} 0 & 0 & 0 & \dfrac{dN_1}{dx} & 0 & 0 & 0 & \dfrac{dN_2}{dx} \\ 0 & 0 & 0 & 0 & 0 & 0 & 0 & 0 \end{bmatrix}, \quad [B_s] = \begin{bmatrix} 0 & 0 & 0 & 0 & 0 & 0 & 0 & 0 \\ 0 & 0 & 0 & N_1 & 0 & 0 & 0 & N_2 \end{bmatrix}$$

By introducing the strain-displacement relation (Eq. 30), the Eq. (22) can be rewritten as:

$$\begin{aligned} 0 = \delta d^T \int_0^L \Big( & B^{0^T} A_{11} B^0 + B^{0^T} B_{11} B^1 + B^{0^T} B_{11}^s B^2 + B^{1^T} B_{11} B^0 + \\ & B^{1^T} D_{11} B^1 + B^{1^T} D_{11}^s B^2 + B^{2^T} B_{11}^s B^0 + B^{2^T} D_{11}^s B^1 + \\ & B^{2^T} H_{11}^s B^2 + B^{s^T} A_{55}^s B^s \Big) dx \, \delta d - q \int_0^L \delta w \, dx \end{aligned} \quad (31)$$

After integration and assembly, the equilibrium equation can be expressed as:

$$[K]\{d\} = \{F\} \quad (32)$$

where $\{F\}$ is the element force vector, $[K]$ represents the element elastic stiffness matrix of new beam element:

$$\begin{aligned} [K] = \sum_e \int_0^L \Big( & B^{0^T} A_{11} B^0 + B^{0^T} B_{11} B^1 + B^{0^T} B_{11}^s B^2 + \\ & B^{1^T} B_{11} B^0 + B^{1^T} D_{11} B^1 + B^{1^T} D_{11}^s B^2 + \\ & B^{2^T} B_{11}^s B^0 + B^{2^T} D_{11}^s B^1 + B^{2^T} H_{11}^s B^2 + B^{s^T} A_{55}^s B^s \Big) dx \end{aligned} \quad (33)$$

## 5. Numerical results

In this section, the accuracy of the present FEM solutions is first demonstrated for the FG straight beams to prove its validity; and then extended for the analysis of curved FG beams. In addition, the effects of the power law index $p$, the thickness-to-side ratio and the material properties on the bending behavior of the isotropic FG beam and FG sandwich beam have been investigated.

Typical values for metal and ceramics used in the FG sandwich beam are listed in Table 1. The applied boundary conditions (BCs) considered in the present study are illustrated in Table 2.

## 5.1 Convergence and validation study

At first example, the convergence study of the developed element model is carried-out. In order to assess the validity of the proposed element model, it is necessary to apply it for the straight FG beams and then extended for the curved FG beams. The numerical results for the straight FG beams can be easily obtained by setting $R = \infty$. In this analysis, three common types (A, B and C) of FG beams, described previously, are considered. This example has been investigated by Li et al. [47] by considering simply supported (SS) beam, doubly clamped beam (CC) under distributed load $q_z$ and cantilever beam (CF) under concentrated load $F_z$. Fig. 4 shows the FG beams analyzed, where beam geometry, material distribution, boundary conditions and loads are described. The computed results are obtained for five values of volume fraction index ($p = 0, 0.5, 1, 5, 10$). The convergence of the maximal deflection for SS, CC and CF FG beams is presented, respectively, in Tables 3, 4 and 5 with different mesh sizes (ne = 2, 4, 8, 12, 16, 24, 32). The obtained results are compared to benchmark solutions of Koutoati et al. [54] and Li et al. [55] which uses a finite element models based on the HSDT. It is clear that for all types of FG beams, volume fraction index ($p$) and boundary conditions that the presented results are in excellent agreement with the solutions in existing literature. Thus, the performance of present finite element formulation is ascertained. It is noted that, for SS and CF FG sandwich beams (type B), it needs just two (02) elements to achieve a desirable level of accuracy with the reference solution and for CC FG sandwich beams it needs (08) eight elements for the same.

The effect of volume fraction exponent ($p$) on the transverse displacement of several types of FG beams using various boundary conditions (SS, CC and CF) is discussed. The exponent is chosen as $p = 0, 0.5, 1, 5, 10$. The differences between the results of various boundary conditions are very significant. Also, it can be observed that increasing the value of volume fraction exponents ($p$) increases the center displacement in all sequences and boundary conditions. These results are expected because the larger volume fraction index ($p$) means the beam has a smaller ceramic component whose Young's modulus is greater than that of metal and hence the stiffness is reduced. In addition, it is evident that the maximal deflection decreases as the rigidity of boundary restraint is increased.

## 5.2 Single-layer FG beam (Type A)

In the second example, the validation of the proposed FE beam model is carried out by comparing the obtained results with those computed via the analytical solution based on the HSDT developed by Li et al. [31]; Navier solutions and finite element model developed by Vo et al. [48]. A simply supported single layer FG beam (Type A) subjected to uniformly distributed load is considered. The top surface of the FG beam is ceramic-rich and the bottom surface is metal-rich. The study is performed for different volume fraction index ($p = 0, 1, 2, 5, 10$) and two length-to-height ratio ($L/h = 5$ and 20). For easiness, the following non-dimensional terms are used:

- Non-dimensional transverse displacement:

$$\bar{w} = \begin{cases} \dfrac{100 E_m h^3}{qL^4} w\left(\dfrac{L}{2}, z\right) & \text{for SS and CC beams} \\ \dfrac{100 E_m h^3}{qL^4} w(L, z) & \text{for CF beams} \end{cases} \quad (34)$$

- Non-dimensional axial and shear stresses:

$$\begin{aligned} \bar{\sigma}_x &= \dfrac{h}{qL} \sigma_x\left(\dfrac{L}{2}, z\right) \\ \bar{\tau}_{xz} &= \dfrac{h}{qL} \tau_{xz}(0, z) \end{aligned} \quad (35)$$

Based on the convergence study, it is evident that sixteen elements are more than enough to obtain more accurate results. Therefore, this mesh size is employed for all the problems presented in this work. The non-dimensional results of transverse displacement, axial stress and transverse shear stress predicted by proposed model are summarized in Table 6. It can be seen from the table that the present results are very close to those obtained by Li et al. [31] and Vo et al. [48]. Indeed, for thick FG beam ($L/h = 5$), the maximum percentage error of transverse displacement predicted by developed finite element in comparison with Navier solution of Vo et al. [48] is 0.00006%, 0.0004%, 0.0006%, 0.00005%, 0.00009% with respect to the volume fraction index ($p$) of 0, 1, 2, 5 and 10.

After establishing the performance of the present element model, the static analysis of single-layer FG curved beam is examined for various values of curvature ($R/L$ = 5, 10, 20, 50, 100, ∞). Table 7 and 8 show the values of transverse displacement and stresses, respectively. This example aims to verify the obtained results with Quasi-2D trigonometric solutions of Sayyad and Ghugal [62] considering the "Stretching effect" ($\varepsilon_z \neq 0$). It may be observed that the results of the developed element are in good agreement with those reported by Sayyad and Ghugal [62]. There is a little difference between the results of the present FE beam model and Quasi-2D solutions. This is due to the different approaches used to predict the response of the FG beam. Nevertheless, a good agreement between the results is found.

It can be seen from the Tables 7 and 8 that increasing the volume fraction index ($p$) results in an increase in both transverse displacement and axial stress; and decrease in shear stress whatever the radius of curvature. For specific length-to-height ratio, the radius of curvature has a slight effect on the bending response of single layer FG curved beams. Indeed, the values of transverse displacement, axial stress and transverse shear stress are almost the same compared to case of FG straight beams.

The distributions of axial stress and transverse shear stress along the thickness of FG curved beam are plotted in Fig. 5. From the Fig. 5a, the parabolic distribution of transverse shear stress is observed for homogeneous beam (full ceramic and full metal), whereas an asymmetric variation can be seen for FG curved beam. Fig. 5b shows that the variation of axial stress is linear for homogeneous beam and non-linear for FG beams. The variation of transverse shear stress and axial stress of FG curved beam are strongly influenced by the volume fraction index ($p$).

**5.3 Sandwich straight beams with FG face sheets and isotropic core (Type B)**

This example is performed for symmetric and non-symmetric FG straight sandwich beams (Type B). The face sheets are assumed to be made of FG layers. The top and the bottom face sheets are graded from metal to ceramic (Al/Al$_2$O$_3$) and the core layer is made of pure ceramic (Al$_2$O$_3$). The beam is simply supported at both ends (SS) and subjected to uniformly distributed load. In the current study, four different sequences (1-1-1, 1-2-1, 2-1-1, 2-2-1), two length-to-height ratio ($L/h$ = 5 and 20) and five volume fraction index ($p$ = 0, 1, 2, 5, 10) are considered. Tables 9-11 present the comparison of the non-dimensional transverse displacements, axial stress and transverse shear stress, respectively. The comparison results verify the accuracy of the developed element where one can see clearly, for all schemes and both thin and thick beams, that the present results are in excellent agreement with those presented by Vo et al. [48] using Navier solution based on the third shear deformation theory (TSDT) and Sayyad and Avhad [68] using the hyperbolic shear deformation theory (HSDT). Moreover, the effect of volume fraction index ($p$) on the non-dimensional central deflection for different core-to-face sheets thickness ratio ($h_c/h_f$) is illustrated in Fig. 6. It can be seen that the lowest and highest values of deflection correspond to the (1-2-1) and (2-1-1) FG sandwich beams, respectively. It is also observed from these tables that as the core thickness increases the non-dimensional values of deflection and stresses decrease. This is due to high proportion of ceramic which leads the plate to be more rigid.

**5.4 Sandwich curved beams with FG face sheets and isotropic core (Type B)**

In this example, the bending of FG sandwich curved beams (Type B) is investigated for various radii of curvature ($R/L$ = 5, 10, 20, 50, 100, $\infty$). Two types of symmetric (1-1-1) and non-symmetric (2-2-1) sandwich curved beams are considered. The computed results are obtained for two length-to-height ratio ($L/h$ = 5 and 10) and five volume fraction index ($p$ = 0, 1, 2, 5, 10). The non-dimensional results of vertical displacement and stresses of simply supported (1-1-1) FG sandwich curved beam are presented in Tables 12-14. It is found that the present results are still in good agreement with referential results available in the literature [62, 63] by considering the stretching effect.

Table 15 shows some new results for CC and CF (1-1-1) FG sandwich curved beams. From the Tables 12 and 15, it can be seen that the vertical displacement decreases as the rigidity of boundary restraint is increased. Thus, the boundary conditions have more significant effects on the vertical displacement of FG sandwich curved beams. In addition, the numerical results for bending behavior of non-symmetric (2-2-1) FG sandwich SS curved beams are shown in Table 16. It should be noted that the results of tables 15-16 are presented for the first time which will serve as a benchmark for the future investigation.

Figs. 7 and 8 show, respectively, the variation of non-dimensional stresses ($\bar{\sigma}_{xx}, \bar{\tau}_{xz}$) through the thickness of (1-1-1) and (2-2-1) FG sandwich SS curved beams. It can be seen that there are some differences in the distribution of stresses between the symmetric and non-symmetric FG sandwich curved beams. In Figs. 7a and 8a, it is observed for both symmetric and non-symmetric FG sandwich beams that the maximum shear stress occurs at the middle plane of the beam. Further, as the volume fraction index ($p$) increases, the non-dimensional transverse shear stress values increase. From the Fig. 7b, it is interesting to see for symmetric (1-1-1) beams that the same maximum compressive (tensile) axial stress is located at the top and bottom face sheets, while the variation of axial stress in non-symmetric (2-2-1) beams is not the same as that found in the symmetric beams (Fig. 8b). This difference is caused by the gradation in material properties. In addition, it is observed, for all values of volume fraction index, that the variation of axial stress along the thickness of top and bottom face sheets is nonlinear, and linear through the thickness of the core. Moreover, as the volume fraction index ($p$) increases, the non-dimensional axial stress values decrease.

**5.5 Sandwich beams with FG core and homogeneous face sheets (Type C)**

Finally, the bending response of (1-8-1) sandwich beam of Type C under uniformly distributed load is analyzed. The top and the bottom face sheets are made of pure metal (Al) and ceramic ($Al_2O_3$), respectively, whereas the FG core layer is graded from ceramic to metal ($Al_2O_3$/Al). Different volume fraction index ($p = 0.5, 1, 2, 5, 10$) are considered for the investigation. The results of the non-dimensional transverse displacement with different boundary condition and stresses of straight sandwich beams are displayed in Tables 17 to 19. As expected, the comparison results show the accuracy of the developed element where one can see clearly that the present results are in excellent agreement with referential results of Vo et al. [48].

In Table 19, the non-dimensional results of transverse displacement and stresses of (1-8-1) FG sandwich SS curved beams are given for the first time. The numerical results are obtained for various radii of curvature ($R/L = 5, 10, 20, 50, 100, \infty$), two length-to-height ratio ($L/h = 5$ and 10) and various volume fraction index ($p = 0, 1, 2, 5, 10$). It can be seen from the table that as the radius of curvature ($R/L$) increases, the transverse displacement and axial stress increase slightly up to $R/L \leq 50$ while further increasing this ratio ($R/L > 50$) has no remarkable effect on the transverse displacement and stresses. On the other hand, the values of transverse shear stresses are almost the same compared to case of FG straight beams.

In Fig. 9, the variation of the non-dimensional axial stress and transverse shear stress through the thickness is plotted for various volume fraction index ($p$). It can be seen from this type of FG sandwich beams (Type C) with $p = 10$ that the maximum shear stress is obtained at the top surface of the core layer (Fig. 9a), while the maximum axial stress is observed around the top face sheet (ceramic rich) (see Fig. 9b).

## 6. Conclusion

In this paper, a new higher-order shear deformation theory is proposed to study the static behavior of FG curved sandwich beam. The present theory is proven to provide an accurate distribution of transverse shear stress through the thickness of FG beam and satisfies the zero traction boundary conditions on the top and bottom surfaces without using any shear correction factors. Based on the suggested model, a new efficient two-noded FG beam element with 4 DOFs is successfully developed for the first time to determine accurately the displacement and stresses of FG curved sandwich beam. Three common configurations of FG beams are used for the study, namely: (a) single layer FG beam; (b) sandwich beam with FG face sheets and homogeneous core and (c) sandwich beams with homogeneous face sheets and FG core. The performance and reliability of the developed finite element model for analyzing the bending behavior of FG curved sandwich beam are validated with existing literature. The advantage of present element model is seen to enable the solution of a variety of problems considering symmetric and non-symmetric FG curved sandwich beams with various boundary conditions and arbitrary material distribution. Effects of the volume fraction index, radius of curvature, material distributions, length-to-thickness ratio, face-to-core- thickness ratio, loadings and boundary conditions on the deflection and stresses are discussed. The obtained results are compared with analytical solutions and those predicted by state of the art advanced finite element models available in the literature. The comparison shows the accuracy and fast rate of convergence of the proposed finite element model. Further, it can be deduced that the proposed model is able to predict accurately the deflection and stresses of thin and thick straight FG sandwich beams as well as those of curved sandwich beams. The important key points that can be concluded from this investigation are summarized as follows:

- For all types of boundary conditions, the volume fraction index ($p$) significantly affect the deflection and stresses of the FG sandwich curved beams. As the volume fraction index increases, the non-dimensional displacements and axial stress increase and transverse shear stress decreases significantly whatever the radius of curvature.

- For specific length-to-height ratio, the radius of curvature has a slight effect on the bending analysis of FG sandwich curved beams.

- With increase in radius of curvature of beam, the values of transverse displacement, axial stress and transverse shear stresses are almost the same compared to case of FG straight beams.

- The core thickness has a significant effect on the mechanical properties of FG sandwich curved beams compared to that of the face layers. As the core thickness increases, the deflection value decreases. This is due to high proportion of ceramic which leads the beam to be more rigid.

Finally, it can be concluded that the proposed model is not only accurate and efficient but also simplifies the application in predicting the static response of FG sandwich curved beams.

**Acknowledgment**

This research was supported by the Algerian Directorate General of Scientific Research and Technological Development (DGRSDT), in Algeria.

**Table Captions**

**Figure captions**

**Fig.1** Geometry and coordinate of isotropic and FG sandwich curved beams.

**Fig. 2** Variation of the volume fraction function through the thickness of three types of FG~~M~~ beams for various values of the power law index (*p*).

**Fig. 3** Present two-noded beam element with corresponding DOFs.

**Fig.4** Studied FG beams with different boundary conditions, material distribution and loads.

**Fig.5** Distribution of non-dimensional stresses through the thickness of single layer FG SS curved beam (Type A), (a) axial stress, (b) transverse shear stress.

**Fig.6** Effect of volume fraction index (*p*) with different core-to-face sheets thickness ratio on the non-dimensional center deflection of SS sandwich straight beams with FG face sheets ($a/h = 5$).

**Fig.7** Distribution of non-dimensional stresses along the thickness of symmetric (1-1-1) FG sandwich SS curved beams (Type B), (a) axial stress, (b) transverse shear stress

**Fig.8** Distribution of non-dimensional stresses along the thickness of non-symmetric (2-2-1) FG sandwich SS curved beams (Type B), (a) axial stress, (b) transverse shear stress.

**Fig.9** Distribution of non-dimensional stresses along the thickness of (1-8-1) FG sandwich SS curved beams (Type C), (a) axial stress, (b) transverse shear stress.

**Table1.** Material properties used in the FG sandwich beam.

| Properties | Metal: Ti–6A1–4V | Ceramic: ZrO$_2$ |
|---|---|---|
| $E$ (GPa) | 70.0 | 380.0 |
| $\nu$ | 0.3 | |

**Table 2.** Boundary conditions used in the present study.

| Boundary conditions | Left boundary ($x = 0$) | Right boundary ($x = L$) |
|---|---|---|
| Simply supported (SS) | $w_0 = 0,\ u_0 \neq 0, w_{0,x} \neq 0, \varphi_x \neq 0$ | $w_0 = 0,\ u_0 \neq 0, w_{0,x} \neq 0, \varphi_x \neq 0$ |
| Clamped-Clamped (CC) | $w_0 = 0,\ u_0 = w_{0,x} = \varphi_x = 0$ | $w_0 = 0,\ u_0 = w_{0,x} = \varphi_x = 0$ |
| Clamped-Free (CF) | $w_0 = 0,\ u_0 = w_{0,x} = \varphi_x = 0$ | $w_0 \neq 0,\ u_0 \neq 0, w_{0,x} \neq 0, \varphi_x \neq 0$ |

**Table 3** Beam theories comparisons of deflection at L/2 of the double simply supported beam (SS) in terms of the volume fraction index (*p*) (values in mm).

| Type | References | Model | Volume fraction index (*p*) | | | | |
|---|---|---|---|---|---|---|---|
| | | | 0 | 0.5 | 1 | 5 | 10 |
| A (Full FGM) | Present (ne=2) | FE-PSDT | 84.287 | 127.90 | 163.23 | 247.17 | 276.94 |
| | Present (ne=4) | FE-PSDT | 84.288 | 129.24 | 167.14 | 255.11 | 282.26 |
| | Present (ne=8) | FE-PSDT | 84.288 | 129.57 | 168.12 | 257.10 | 283.58 |
| | Present (ne=12) | FE-PSDT | 84.288 | 129.63 | 168.30 | 257.44 | 283.83 |
| | Present (ne=16) | FE-PSDT | 84.288 | 129.66 | 168.37 | 257.60 | 283.92 |
| | Present (ne=24) | FE-PSDT | 84.288 | 129.67 | 168.41 | 257.69 | 283.98 |
| | Present (ne=32) | FE-PSDT | 84.288 | 129.68 | 168.43 | 257.72 | 284.00 |
| | Koutoati et al. [52] | HSDT | 84.290 | 129.64 | 168.45 | 257.72 | 284.01 |
| | Li et al. [53] | DTS | 84.289 | 129.63 | 168.45 | 257.73 | 284.01 |
| B (3-4-3) | Present (ne=2) | FE-PSDT | 84.288 | 126.61 | 162.00 | 281.10 | 314.71 |
| | Present (ne=4) | FE-PSDT | 84.288 | 126.61 | 162.00 | 281.10 | 314.71 |
| | Present (ne=8) | FE-PSDT | 84.288 | 126.61 | 162.00 | 281.10 | 314.71 |
| | Present (ne=12) | FE-PSDT | 84.288 | 126.61 | 162.00 | 281.10 | 314.71 |
| | Present (ne=16) | FE-PSDT | 84.288 | 126.61 | 162.00 | 281.10 | 314.71 |
| | Present (ne=24) | FE-PSDT | 84.288 | 126.61 | 162.00 | 281.10 | 314.71 |
| | Present (ne=32) | FE-PSDT | 84.288 | 126.61 | 162.00 | 281.10 | 314.71 |
| | Koutoati et al. [52] | HSDT | 84.290 | 126.60 | 162.00 | 281.12 | 314.74 |
| | Li et al. [53] | DTS | 84.289 | 126.59 | 162.00 | 281.12 | 314.74 |
| C (3-4-3) | Present (ne=2) | FE-PSDT | 159.79 | 182.16 | 197.28 | 213.89 | 213.89 |
| | Present (ne=4) | FE-PSDT | 164.38 | 190.01 | 204.04 | 226.52 | 228.81 |
| | Present (ne=8) | FE-PSDT | 165.53 | 191.98 | 206.48 | 229.68 | 231.99 |
| | Present (ne=12) | FE-PSDT | 165.74 | 192.34 | 206.93 | 230.27 | 232.58 |
| | Present (ne=16) | FE-PSDT | 165.81 | 192.47 | 207.09 | 230.47 | 232.79 |
| | Present (ne=24) | FE-PSDT | 165.87 | 192.56 | 207.20 | 230.62 | 232.94 |
| | Present (ne=32) | FE-PSDT | 165.90 | 192.59 | 207.24 | 230.67 | 232.99 |
| | Koutoati et al. [52] | HSDT | 166.35 | 192.63 | 207.30 | 230.68 | 232.98 |
| | Li et al. [53] | DTS | 166.35 | 192.63 | 207.31 | 230.69 | 232.98 |

**Table 4** Beam theories comparisons of deflection at L/2 of the double clamped beam (CC) in terms of the volume fraction index ($p$) (values in mm).

| Type | References | Model | Volume fraction index ($p$) | | | | |
|---|---|---|---|---|---|---|---|
| | | | 0 | 0.5 | 1 | 5 | 10 |
| A (Full FGM) | Present (ne=2) | FE-PSDT | 16.147 | 23.604 | 27.777 | 39.414 | 47.819 |
| | Present (ne=4) | FE-PSDT | 17.983 | 27.016 | 34.287 | 53.151 | 60.267 |
| | Present (ne=8) | FE-PSDT | 18.201 | 27.645 | 35.634 | 55.958 | 62.606 |
| | Present (ne=12) | FE-PSDT | 18.295 | 27.834 | 35.974 | 56.678 | 63.286 |
| | Present (ne=16) | FE-PSDT | 18.343 | 27.920 | 36.118 | 56.983 | 63.589 |
| | Present (ne=24) | FE-PSDT | 18.389 | 28.000 | 36.240 | 57.243 | 63.857 |
| | Present (ne=32) | FE-PSDT | 18.410 | 28.038 | 36.292 | 57.352 | 63.973 |
| | Koutoati et al. [52] | HSDT | 18.450 | 28.090 | 36.390 | 57.500 | 64.160 |
| | Li et al. [53] | DTS | 18.454 | 28.088 | 36.387 | 57.503 | 64.162 |
| B (3-4-3) | Present (ne=2) | FE-PSDT | 16.447 | 24.865 | 31.910 | 55.658 | 62.363 |
| | Present (ne=4) | FE-PSDT | 17.983 | 26.582 | 33.745 | 57.767 | 64.536 |
| | Present (ne=8) | FE-PSDT | 18.201 | 26.826 | 34.007 | 58.067 | 64.846 |
| | Present (ne=12) | FE-PSDT | 18.296 | 26.931 | 34.120 | 58.198 | 64.980 |
| | Present (ne=16) | FE-PSDT | 18.343 | 26.985 | 34.178 | 58.264 | 65.048 |
| | Present (ne=24) | FE-PSDT | 18.389 | 27.037 | 34.234 | 58.330 | 65.116 |
| | Present (ne=32) | FE-PSDT | 18.410 | 27.062 | 34.261 | 58.361 | 65.148 |
| | Koutoati et al. [52] | HSDT | 18.450 | 27.110 | 34.320 | 58.440 | 65.240 |
| | Li et al. [53] | DTS | 18.454 | 27.109 | 34.320 | 58.443 | 65.240 |
| C (3-4-3) | Present (ne=2) | FE-PSDT | 26.582 | 27.466 | 27.777 | 28.312 | 28.543 |
| | Present (ne=4) | FE-PSDT | 32.966 | 37.521 | 40.050 | 44.646 | 45.403 |
| | Present (ne=8) | FE-PSDT | 34.368 | 39.797 | 42.847 | 48.329 | 49.170 |
| | Present (ne=12) | FE-PSDT | 34.691 | 40.296 | 43.452 | 49.140 | 50.011 |
| | Present (ne=16) | FE-PSDT | 34.821 | 40.491 | 43.688 | 49.457 | 50.343 |
| | Present (ne=24) | FE-PSDT | 34.930 | 40.648 | 43.876 | 49.712 | 50.610 |
| | Present (ne=32) | FE-PSDT | 34.973 | 40.711 | 43.951 | 49.813 | 50.717 |
| | Koutoati et al. [52] | HSDT | 35.150 | 40.820 | 44.080 | 49.930 | 50.820 |
| | Li et al. [53] | DTS | 35.150 | 40.824 | 44.081 | 49.928 | 50.815 |

**Table 5** Beam theories comparisons of deflection at the tip end of the cantilever beam (CF) in terms of the power law index $p$ (values in mm).

| Type | Reference | Model | Volume fraction index ($p$) | | | | |
|---|---|---|---|---|---|---|---|
| | | | 0 | 0.5 | 1 | 5 | 10 |
| A (Full FGM) | Present (ne=2) | FE-PSDT | 13.509 | 20.692 | 26.730 | 40.800 | 45.199 |
| | Present (ne=4) | FE-PSDT | 13.538 | 20.798 | 26.975 | 41.306 | 45.598 |
| | Present (ne=8) | FE-PSDT | 13.552 | 20.834 | 27.048 | 41.459 | 45.730 |
| | Present (ne=12) | FE-PSDT | 13.557 | 20.843 | 27.065 | 41.494 | 45.763 |
| | Present (ne=16) | FE-PSDT | 13.559 | 20.848 | 27.072 | 41.508 | 45.777 |
| | Present (ne=24) | FE-PSDT | 13.561 | 20.851 | 27.077 | 41.520 | 45.789 |
| | Present (ne=32) | FE-PSDT | 13.562 | 20.852 | 27.079 | 41.524 | 45.793 |
| | Koutoati et al. [52] | HSDT | 13.560 | 20.850 | 27.080 | 41.520 | 45.800 |
| | Li et al. [53] | DTS | 13.560 | 20.840 | 27.080 | 41.520 | 45.800 |
| B (3-4-3) | Present (ne=2) | FE-PSDT | 13.509 | 20.285 | 25.948 | 45.008 | 50.387 |
| | Present (ne=4) | FE-PSDT | 13.538 | 20.317 | 25.982 | 45.048 | 50.428 |
| | Present (ne=8) | FE-PSDT | 13.552 | 20.333 | 26.000 | 45.068 | 50.449 |
| | Present (ne=12) | FE-PSDT | 13.557 | 20.338 | 26.005 | 45.075 | 50.455 |
| | Present (ne=16) | FE-PSDT | 13.560 | 20.341 | 26.008 | 45.078 | 50.459 |
| | Present (ne=24) | FE-PSDT | 13.561 | 20.343 | 26.011 | 45.081 | 50.462 |
| | Present (ne=32) | FE-PSDT | 13.562 | 20.344 | 26.012 | 45.082 | 50.463 |
| | Koutoati et al. [52] | HSDT | 13.564 | 20.342 | 26.014 | 45.088 | 50.471 |
| | Li et al. [53] | DTS | 13.564 | 20.343 | 26.015 | 45.088 | 50.470 |
| C (3-4-3) | Present (ne=2) | FE-PSDT | 26.267 | 30.331 | 32.555 | 36.131 | 36.503 |
| | Present (ne=4) | FE-PSDT | 26.530 | 30.765 | 33.090 | 36.832 | 37.217 |
| | Present (ne=8) | FE-PSDT | 26.604 | 30.884 | 33.236 | 37.025 | 37.415 |
| | Present (ne=12) | FE-PSDT | 26.620 | 30.909 | 33.266 | 37.065 | 37.456 |
| | Present (ne=16) | FE-PSDT | 26.627 | 30.918 | 33.277 | 37.080 | 37.472 |
| | Present (ne=24) | FE-PSDT | 26.632 | 30.926 | 33.286 | 37.092 | 37.484 |
| | Present (ne=32) | FE-PSDT | 26.634 | 30.928 | 33.289 | 37.096 | 37.489 |
| | Koutoati et al. [52] | HSDT | 26.707 | 30.933 | 33.296 | 37.094 | 37.481 |
| | Li et al. [53] | DTS | 26.707 | 30.933 | 33.296 | 37.093 | 37.481 |

**Table 6.** Comparison of the maximum transverse displacement, axial stress and shear stress of simply supported FG straight beam (Type A).

| L/h | p | Reference | Model | $\bar{w}(L/2,0)$ | $\bar{\sigma}_{xx}(L/2,h/2)$ | $\bar{\tau}_{xz}(0,0)$ |
|---|---|---|---|---|---|---|
| 5 | 0 | Present | FE-PSDT | 3.1652 | 3.8136 | 0.7534 |
|  |  | Vo et al. [46] | Navier | 3.1654 | 3.8020 | 0.7332 |
|  |  | Li et al. [29] | HSDT | 3.1657 | 3.8020 | 0.7500 |
|  |  | Vo et al. [46] | FE-TBT | 3.1654 | 3.8021 | 0.7335 |
|  | 1 | Present | FE-PSDT | 6.2563 | 5.9061 | 0.7534 |
|  |  | Vo et al. [46] | Navier | 6.2594 | 5.8837 | 0.7332 |
|  |  | Li et al. [29] | HSDT | 6.2599 | 5.8836 | 0.7500 |
|  |  | Vo et al. [46] | FE-TBT | 6.2590 | 5.8870 | 0.7335 |
|  | 2 | Present | FE-PSDT | 8.0628 | 6.9090 | 0.6908 |
|  |  | Vo et al. [46] | Navier | 8.0677 | 6.8826 | 0.6706 |
|  |  | Li et al. [29] | HSDT | 8.0602 | 6.8812 | 0.6787 |
|  |  | Vo et al. [46] | FE-TBT | 8.0668 | 6.8860 | 0.6700 |
|  | 5 | Present | FE-PSDT | 9.8276 | 8.1460 | 0.6111 |
|  |  | Vo et al. [46] | Navier | 9.8281 | 8.1106 | 0.5905 |
|  |  | Li et al. [29] | HSDT | 9.7802 | 8.1030 | 0.5790 |
|  |  | Vo et al. [46] | FE-TBT | 9.8271 | 8.1150 | 0.5907 |
|  | 10 | Present | FE-PSDT | 10.937 | 9.7544 | 0.6675 |
|  |  | Vo et al. [46] | Navier | 10.938 | 9.7122 | 0.6467 |
|  |  | Li et al. [29] | HSDT | 10.897 | 9.7063 | 0.6436 |
|  |  | Vo et al. [46] | FE-TBT | 10.937 | 9.7170 | 0.6477 |
| 20 | 0 | Present | FE-PSDT | 2.8962 | 15.0525 | 0.7626 |
|  |  | Vo et al. [46] | Navier | 2.8962 | 15.0129 | 0.7451 |
|  |  | Li et al. [29] | HSDT | 2.8962 | 15.0130 | 0.7500 |
|  |  | Vo et al. [46] | FE-TBT | 2.8963 | 15.0200 | 0.7470 |
|  | 1 | Present | FE-PSDT | 5.8021 | 23.2833 | 0.7626 |
|  |  | Vo et al. [46] | Navier | 5.8049 | 23.2053 | 0.7451 |
|  |  | Li et al. [29] | HSDT | 5.8049 | 23.2054 | 0.7500 |
|  |  | Vo et al. [46] | FE-TBT | 5.8045 | 23.2200 | 0.7470 |
|  | 2 | Present | FE-PSDT | 7.4366 | 27.1888 | 0.7003 |
|  |  | Vo et al. [46] | Navier | 7.4421 | 27.0991 | 0.6824 |
|  |  | Li et al. [29] | HSDT | 7.4415 | 27.0989 | 0.6787 |
|  |  | Vo et al. [46] | FE-TBT | 7.4412 | 27.1100 | 0.6777 |
|  | 5 | Present | FE-PSDT | 8.8128 | 31.9301 | 0.6212 |
|  |  | Vo et al. [46] | Navier | 8.8182 | 31.8130 | 0.6023 |
|  |  | Li et al. [29] | HSDT | 8.8151 | 31.8112 | 0.5790 |
|  |  | Vo et al. [46] | FE-TBT | 8.8173 | 31.8300 | 0.6039 |
|  | 10 | Present | FE-PSDT | 9.6868 | 38.2826 | 0.6785 |
|  |  | Vo et al. [46] | Navier | 9.6905 | 38.1385 | 0.6596 |
|  |  | Li et al. [29] | HSDT | 9.6879 | 38.1372 | 0.6436 |
|  |  | Vo et al. [46] | FE-TBT | 9.6899 | 38.1600 | 0.6682 |

**Table 7.** Non-dimensional maximum deflection of simply supported FG curved beam (Type A).

| $L/h$ | $p$ | Model | $\varepsilon_z$ | R/L | | | | | |
|---|---|---|---|---|---|---|---|---|---|
| | | | | 5 | 10 | 20 | 50 | 100 | $\infty$ |
| 5 | 0 | Present | $=0$ | 3.1638 | 3.1649 | 3.1651 | 3.1652 | 3.1652 | 3.1652 |
| | | SSDT [60] | $\neq 0$ | 3.1355 | 3.1356 | 3.1357 | 3.1357 | 3.1357 | 3.1357 |
| | 1 | Present | $=0$ | 6.2480 | 6.2528 | 6.2547 | 6.2557 | 6.2560 | 6.2563 |
| | | SSDT [60] | $\neq 0$ | 6.1482 | 6.1316 | 6.1233 | 6.1184 | 6.1168 | 6.1151 |
| | 2 | Present | $=0$ | 8.0513 | 8.0578 | 8.0605 | 8.0619 | 8.0624 | 8.0628 |
| | | SSDT [60] | $\neq 0$ | 7.8773 | 7.8558 | 7.8453 | 7.8389 | 7.8368 | 7.8347 |
| | 5 | Present | $=0$ | 9.8153 | 9.8223 | 9.8251 | 9.8266 | 9.8271 | 9.8276 |
| | | SSDT [60] | $\neq 0$ | 9.6298 | 9.6069 | 9.5956 | 9.5889 | 9.5866 | 9.5844 |
| | 10 | Present | $=0$ | 10.925 | 10.932 | 10.934 | 10.936 | 10.936 | 10.937 |
| | | SSDT [60] | $\neq 0$ | 10.787 | 10.761 | 10.748 | 10.740 | 10.738 | 10.735 |
| 10 | 0 | Present | $=0$ | 2.9453 | 2.9489 | 2.9498 | 2.9500 | 2.9501 | 2.9501 |
| | | SSDT [60] | $\neq 0$ | 2.9396 | 2.9397 | 2.9398 | 2.9398 | 2.9398 | 2.9398 |
| | 1 | Present | $=0$ | 5.8729 | 5.8853 | 5.8897 | 2.9500 | 5.8924 | 5.8930 |
| | | SSDT [60] | $\neq 0$ | 5.8188 | 5.8017 | 5.7898 | 5.7881 | 5.7864 | 5.7847 |
| | 2 | Present | $=0$ | 7.5354 | 7.5514 | 7.5573 | 7.5603 | 7.5611 | 7.5620 |
| | | SSDT [60] | $\neq 0$ | 7.4144 | 7.3884 | 7.3815 | 7.3750 | 7.3729 | 7.3697 |
| | 5 | Present | $=0$ | 8.9883 | 9.0050 | 9.0112 | 9.0142 | 9.0151 | 9.0160 |
| | | SSDT [60] | $\neq 0$ | 8.8619 | 8.8386 | 8.8270 | 8.8200 | 8.8177 | 8.8154 |
| | 10 | Present | $=0$ | 9.9110 | 9.9271 | 9.9329 | 9.9357 | 9.9364 | 9.9372 |
| | | SSDT [60] | $\neq 0$ | 9.8430 | 9.8163 | 9.8030 | 9.7951 | 9.7924 | 9.7898 |

**Table 8.** Non-dimensional transverse shear and axial stresses of simply supported FG curved beam (Type A).

| L/h | p | Model | $\varepsilon_z$ | R/L | | | | | |
|---|---|---|---|---|---|---|---|---|---|
| | | | | 5 | 10 | 20 | 50 | 100 | ∞ |
| Shear stress ($\bar{\tau}_{xz}$) | | | | | | | | | |
| 5 | 0 | Present | = 0 | 0.7532 | 0.7534 | 0.7534 | 0.7534 | 0.7534 | 0.7534 |
| | | SSDT [60] | ≠ 0 | 0.7436 | 0.7436 | 0.7436 | 0.7436 | 0.7436 | 0.7436 |
| | 1 | Present | = 0 | 0.7528 | 0.7532 | 0.7532 | 0.7534 | 0.7534 | 0.7534 |
| | | SSDT [60] | ≠ 0 | 0.7432 | 0.7431 | 0.7430 | 0.7430 | 0.7430 | 0.7430 |
| | 2 | Present | = 0 | 0.6902 | 0.6906 | 0.6906 | 0.6907 | 0.6908 | 0.6908 |
| | | SSDT [60] | ≠ 0 | 0.6832 | 0.6955 | 0.6990 | 0.6997 | 0.6998 | 0.6830 |
| | 5 | Present | = 0 | 0.6107 | 0.6109 | 0.6109 | 0.6111 | 0.6111 | 0.6111 |
| | | SSDT [60] | ≠ 0 | 0.6074 | 0.6073 | 0.6073 | 0.6072 | 0.6072 | 0.6072 |
| | 10 | Present | = 0 | 0.6671 | 0.6673 | 0.6673 | 0.6675 | 0.6675 | 0.6675 |
| | | SSDT [60] | ≠ 0 | 0.6619 | 0.6619 | 0.6618 | 0.6618 | 0.6618 | 0.6618 |
| 10 | 0 | Present | = 0 | 0.7595 | 0.7602 | 0.7604 | 0.7605 | 0.7605 | 0.7605 |
| | | SSDT [60] | ≠ 0 | 0.7566 | 0.7566 | 0.7566 | 0.7566 | 0.7566 | 0.7566 |
| | 1 | Present | = 0 | 0.7588 | 0.7599 | 0.7603 | 0.7604 | 0.7605 | 0.7605 |
| | | SSDT [60] | ≠ 0 | 0.7562 | 0.7561 | 0.7560 | 0.7560 | 0.7560 | 0.7560 |
| | 2 | Present | = 0 | 0.6966 | 0.6975 | 0.6979 | 0.6980 | 0.6981 | 0.6981 |
| | | SSDT [60] | ≠ 0 | 0.6955 | 0.6954 | 0.6954 | 0.6953 | 0.6953 | 0.6953 |
| | 5 | Present | = 0 | 0.6176 | 0.6184 | 0.6186 | 0.6187 | 0.6188 | 0.6188 |
| | | SSDT [60] | ≠ 0 | 0.6189 | 0.6188 | 0.6188 | 0.6188 | 0.6187 | 0..6187 |
| | 10 | Present | = 0 | 0.6748 | 0.6755 | 0.6758 | 0.6759 | 0.6759 | 0.6759 |
| | | SSDT [60] | ≠ 0 | 0.6745 | 0.6745 | 0.6744 | 0.6744 | 0.6744 | 0.6744 |
| Axial stress ($\bar{\sigma}_{xx}$) | | | | | | | | | |
| 5 | 0 | Present | = 0 | 3.8172 | 3.8158 | 3.8148 | 3.8141 | 3.8139 | 3.8136 |
| | | SSDT [60] | ≠ 0 | 3.8220 | 3.8401 | 3.8491 | 3.8545 | 3.8563 | 3.8581 |
| | 1 | Present | = 0 | 5.9005 | 5.9040 | 5.9040 | 5.9058 | 5.9060 | 5.9061 |
| | | SSDT [60] | ≠ 0 | 5.9366 | 5.9583 | 5.9689 | 5.9750 | 5.9771 | 5.9791 |
| | 2 | Present | = 0 | 6.9000 | 6.9052 | 6.9052 | 6.9083 | 6.9087 | 6.9090 |
| | | SSDT [60] | ≠ 0 | 6.9448 | 6.9679 | 6.9791 | 6.9856 | 6.9877 | 6.9899 |
| | 5 | Present | = 0 | 8.1387 | 8.1431 | 8.1431 | 8.1455 | 8.1458 | 8.1460 |
| | | SSDT [60] | ≠ 0 | 8.1850 | 8.2117 | 8.2246 | 8.2321 | 8.2346 | 8.2371 |
| | 10 | Present | = 0 | 9.7511 | 9.7536 | 9.7536 | 9.7544 | 9.7544 | 9.7544 |
| | | SSDT [60] | ≠ 0 | 9.7885 | 9.8207 | 9.8362 | 9.8454 | 9.8485 | 9.8515 |
| 10 | 0 | Present | = 0 | 7.5538 | 7.5533 | 7.5507 | 7.5483 | 7.5474 | 7.5465 |
| | | SSDT [60] | ≠ 0 | 7.6056 | 7.6056 | 7.6253 | 7.6371 | 7.6411 | 7.6450 |
| | 1 | Present | = 0 | 11.644 | 11.665 | 11.671 | 11.674 | 11.675 | 11.675 |
| | | SSDT [60] | ≠ 0 | 11.745 | 11.792 | 11.815 | 11.829 | 11.833 | 11.837 |
| | 2 | Present | = 0 | 13.594 | 13.621 | 13.631 | 13.636 | 13.637 | 13.639 |
| | | SSDT [60] | ≠ 0 | 13.716 | 13.766 | 13.791 | 13.805 | 13.810 | 13.814 |
| | 5 | Present | = 0 | 15.992 | 16.016 | 16.024 | 16.028 | 16.029 | 16.030 |
| | | SSDT [60] | ≠ 0 | 16.105 | 16.163 | 16.191 | 16.208 | 16.213 | 16.219 |
| | 10 | Present | = 0 | 19.189 | 19.208 | 19.213 | 19.214 | 19.214 | 19.214 |
| | | SSDT [60] | ≠ 0 | 19.279 | 19.349 | 19.383 | 19.403 | 19.410 | 19.417 |

**Table 9.** Non-dimensional vertical displacement of FG sandwich SS straight beams (Type B).

| p | Model | L/h = 5 | | | | L/h = 20 | | | |
|---|---|---|---|---|---|---|---|---|---|
| | | 1-1-1 | 1-2-1 | 2-1-1 | 2-2-1 | 1-1-1 | 1-2-1 | 2-1-1 | 2-2-1 |
| 0 | Present | 3.1652 | 3.1652 | 3.1652 | 3.1652 | 2.8962 | 2.8962 | 2.8962 | 2.8962 |
| | Navier [46] | 3.1654 | 3.1654 | 3.1654 | 3.1657 | 2.8963 | 2.8963 | 2.8963 | 2.8947 |
| | HSDT [66] | 3.1241 | 3.1241 | - | - | 2.8585 | 2.8585 | - | - |
| 1 | Present | 6.2688 | 5.4125 | 6.5440 | 5.8399 | 5.9400 | 5.1006 | 6.1973 | 5.5158 |
| | Navier [46] | 6.2693 | 5.4122 | 6.5450 | 5.8403 | 5.9401 | 5.1006 | 6.1977 | 5.5161 |
| | HSDT [66] | 6.3011 | 5.0341 | - | - | 5.9561 | 5.3415 | - | - |
| 2 | Present | 8.3880 | 6.7581 | 8.8871 | 7.5570 | 8.0312 | 6.4276 | 8.4991 | 7.2072 |
| | Navier [46] | 8.3893 | 6.7579 | 8.8896 | 7.5583 | 8.0313 | 6.4276 | 8.5000 | 7.2080 |
| | HSDT [66] | 8.2734 | 6.3359 | - | - | 7.9201 | 6.6697 | - | - |
| 5 | Present | 11.2242 | 8.5134 | 11.8189 | 9.7885 | 10.8374 | 8.1642 | 11.3756 | 9.4103 |
| | Navier [46] | 11.2274 | 8.5137 | 11.8246 | 9.7919 | 10.8376 | 8.1642 | 11.3782 | 9.4120 |
| | HSDT [66] | 11.0708 | 8.0576 | - | - | 10.6766 | 8.4045 | - | - |
| 10 | Present | 12.5612 | 9.4041 | 13.0064 | 10.8439 | 12.1590 | 9.0470 | 12.5249 | 10.4503 |
| | Navier [46] | 12.5659 | 9.4050 | 13.0135 | 10.8486 | 12.1593 | 9.0471 | 12.5281 | 10.4526 |
| | HSDT [66] | 12.3910 | 8.9290 | - | - | 12.1030 | 9.2824 | - | - |

**Table 10.** Non-dimensional axial stress of FG sandwich SS straight beams (Type B).

| $p$ | Model | $L/h = 5$ | | | | $L/h = 20$ | | | |
|---|---|---|---|---|---|---|---|---|---|
| | | 1-1-1 | 1-2-1 | 2-1-1 | 2-2-1 | 1-1-1 | 1-2-1 | 2-1-1 | 2-2-1 |
| 0 | Present | 3.8136 | 3.8136 | 3.8136 | 3.8136 | 15.0525 | 15.0525 | 15.0525 | 15.0525 |
| | Navier [46] | 3.8020 | 3.8020 | 3.8020 | 3.8020 | 15.0129 | 15.0129 | 15.0129 | 15.0129 |
| | HSDT [66] | 3.8025 | 3.8025 | - | - | 15.0136 | 15.0136 | - | - |
| 1 | Present | 1.4391 | 1.2366 | 1.3931 | 1.2517 | 5.6999 | 4.8929 | 5.5128 | 4.9513 |
| | Navier [46] | 1.4349 | 1.2329 | 1.3884 | 1.2474 | 5.6850 | 4.8801 | 5.4960 | 4.9364 |
| | HSDT [66] | 1.4614 | 1.2331 | - | - | 5.7370 | 4.8802 | - | - |
| 2 | Present | 1.9438 | 1.5574 | 1.8537 | 1.5928 | 7.7114 | 6.1694 | 7.3452 | 6.3082 |
| | Navier [46] | 1.9382 | 1.5527 | 1.8475 | 1.5873 | 7.6912 | 6.1532 | 7.3227 | 6.2889 |
| | HSDT [66] | 1.9369 | 1.5530 | - | - | 7.6154 | 6.1534 | - | - |
| 5 | Present | 2.6197 | 1.9763 | 2.4148 | 2.0261 | 10.4107 | 7.8400 | 9.5795 | 8.0350 |
| | Navier [46] | 2.6123 | 1.9750 | 2.4069 | 2.0194 | 10.3835 | 7.8194 | 9.5508 | 8.0109 |
| | HSDT [66] | 2.6101 | 1.9707 | - | - | 10.2712 | 7.8196 | - | - |
| 10 | Present | 2.9375 | 2.1889 | 2.6381 | 2.2271 | 11.6818 | 8.6893 | 10.4669 | 8.8367 |
| | Navier [46] | 2.9293 | 2.1826 | 2.6296 | 2.2199 | 11.6513 | 8.6665 | 10.4357 | 8.8104 |
| | HSDT [66] | 2.9268 | 2.1829 | - | - | 11.5237 | 8.6667 | - | - |

**Table 11.** Non-dimensional transverse shear stress of FG sandwich SS straight beams (Type B).

| $p$ | Model | $L/h = 5$ | | | | $L/h = 20$ | | | |
|---|---|---|---|---|---|---|---|---|---|
| | | 1-1-1 | 1-2-1 | 2-1-1 | 2-2-1 | 1-1-1 | 1-2-1 | 2-1-1 | 2-2-1 |
| 0 | Present | 0.7534 | 0.7534 | 0.7534 | 0.7534 | 0.7626 | 0.7626 | 0.7626 | 0.7626 |
| | Navier [46] | 0.7332 | 0.7332 | 0.7332 | 0.7332 | 0.7451 | 0.7451 | 0.7451 | 0.7451 |
| | HSDT [66] | 0.7285 | 0.7285 | - | - | 0.7355 | 0.7355 | - | - |
| 1 | Present | 0.8782 | 0.8329 | 0.9299 | 0.8684 | 0.8844 | 0.8390 | 0.9365 | 0.8747 |
| | Navier [46] | 0.8586 | 0.8123 | 0.9088 | 0.8479 | 0.8681 | 0.8215 | 0.9166 | 0.8552 |
| | HSDT [66] | 0.8767 | 0.8056 | - | - | 0.8726 | 0.8106 | - | - |
| 2 | Present | 0.9433 | 0.8696 | 1.0344 | 0.9274 | 0.9491 | 0.8751 | 1.0411 | 0.9333 |
| | Navier [46] | 0.9249 | 0.8493 | 1.0136 | 0.9075 | 0.9344 | 0.8581 | 1.0242 | 0.9168 |
| | HSDT [66] | 0.9170 | 0.8424 | - | - | 0.9222 | 0.8486 | - | - |
| 5 | Present | 1.0280 | 0.9116 | 1.1948 | 1.0047 | 1.0343 | 0.9171 | 1.2022 | 1.0107 |
| | Navier [46] | 1.0125 | 0.8925 | 1.1742 | 0.9859 | 1.0227 | 0.9014 | 1.1862 | 0.9957 |
| | HSDT [66] | 1.0048 | 0.8851 | - | - | 1.0101 | 0.8897 | - | - |
| 10 | Present | 1.0800 | 0.9333 | 1.3090 | 1.0517 | 1.0865 | 0.9390 | 1.3172 | 1.0578 |
| | Navier [46] | 1.0665 | 0.9151 | 1.2875 | 1.0335 | 1.0773 | 0.9243 | 1.3008 | 1.0436 |
| | HSDT [66] | 1.0586 | 0.9083 | - | - | 1.0642 | 0.9128 | - | - |

**Table 12.** Non-dimensional vertical displacement of (1-1-1) FG sandwich SS curved beams (Type B).

| $L/h$ | $p$ | Model | $\varepsilon_z$ | R/L | | | | | |
|---|---|---|---|---|---|---|---|---|---|
| | | | | 5 | 10 | 20 | 50 | 100 | ∞ |
| 5 | 0 | Present | $= 0$ | 3.1595 | 3.1628 | 3.1650 | 3.1652 | 3.1652 | 3.1652 |
| | | SSDT [60] | $\neq 0$ | 3.1294 | 3.1295 | 3.1295 | 3.1296 | 3.1296 | 3.1296 |
| | | HSDT[53] | $\neq 0$ | 3.1775 | 3.1775 | 3.1775 | - | - | - |
| | 1 | Present | $= 0$ | 6.2529 | 6.2641 | 6.2681 | 6.2686 | 6.2687 | 6.2687 |
| | | SSDT [60] | $\neq 0$ | 6.1913 | 6.1916 | 6.1916 | 9.1917 | 6.1917 | 6.1917 |
| | | HSDT[53] | $\neq 0$ | 6.2763 | 6.2763 | 6.2763 | - | - | - |
| | 2 | Present | $= 0$ | 8.3633 | 8.3809 | 8.3870 | 8.3878 | 8.3879 | 8.3880 |
| | | SSDT [60] | $\neq 0$ | 8.2823 | 8.2827 | 8.2828 | 8.2828 | 8.2828 | 8.2828 |
| | | HSDT[53] | $\neq 0$ | 8.3955 | 8.3955 | 8.3955 | - | - | - |
| | 5 | Present | $= 0$ | 11.185 | 11.214 | 11.222 | 11.223 | 11.224 | 11.224 |
| | | SSDT [60] | $\neq 0$ | 11.077 | 11.078 | 11.078 | 11.078 | 11.078 | 11.078 |
| | | HSDT[53] | $\neq 0$ | 11.248 | 11.248 | 11.248 | - | - | - |
| | 10 | Present | $= 0$ | 12.517 | 12.549 | 12.559 | 12.560 | 12.561 | 12.561 |
| | | SSDT [60] | $\neq 0$ | 12.396 | 12.397 | 12.397 | 12.397 | 12.397 | 12.397 |
| | | HSDT[53] | $\neq 0$ | 12.625 | 12.625 | 12.625 | - | - | - |
| 10 | 0 | Present | $= 0$ | 2.9312 | 2.9379 | 2.9470 | 2.9499 | 2.9500 | 2.9501 |
| | | SSDT [60] | $\neq 0$ | 2.9337 | 2.9339 | 2.9339 | 2.9339 | 2.9339 | 2.9339 |
| | 1 | Present | $= 0$ | 5.9486 | 5.9690 | 5.9965 | 6.0054 | 6.0057 | 6.0057 |
| | | SSDT [60] | $\neq 0$ | 5.9706 | 5.9709 | 5.9709 | 5.9709 | 5.9709 | 5.9709 |
| | 2 | Present | $= 0$ | 8.0116 | 8.0440 | 8.0878 | 6.0054 | 8.1024 | 8.1025 |
| | | SSDT [60] | $\neq 0$ | 8.0550 | 8.0554 | 8.0555 | 8.0556 | 8.0556 | 8.0556 |
| | 5 | Present | $= 0$ | 10.773 | 10.891 | 10.891 | 10.913 | 10.914 | 10.914 |
| | | SSDT [60] | $\neq 0$ | 10.847 | 10.847 | 10.847 | 10.847 | 10.847 | 10.847 |
| | 10 | Present | $= 0$ | 12.075 | 12.133 | 12.212 | 12.238 | 12.239 | 12.239 |
| | | SSDT [60] | $\neq 0$ | 12.162 | 12.162 | 12.163 | 12.163 | 12.163 | 12.163 |

**Table 13.** Non-dimensional transverse shear stress of (1-1-1) FG sandwich SS curved beams (Type B).

| L/h | p | Model | $\varepsilon_z$ | R/L |  |  |  |  |  |
|---|---|---|---|---|---|---|---|---|---|
|  |  |  |  | 5 | 10 | 20 | 50 | 100 | ∞ |
| 5 | 0 | Present | = 0 | 0.7566 | 0.7574 | 0.7576 | 0.7576 | 0.7577 | 0.7577 |
|  |  | SSDT [60] | ≠ 0 | 0.7431 | 0.7431 | 0.7431 | 0.7431 | 0.7431 | 0.7431 |
|  |  | HSDT[53] | ≠ 0 | 0.7318 | 0.7318 | 0.7318 | - | - | - |
|  | 1 | Present | = 0 | 0.8794 | 0.8808 | 0.8811 | 0.8812 | 0.8812 | 0.8812 |
|  |  | SSDT [60] | ≠ 0 | 0.8623 | 0.8623 | 0.8623 | 0.8623 | 0.8623 | 0.8623 |
|  |  | HSDT[53] | ≠ 0 | 0.9090 | 0.9090 | 0.9090 | - | - | - |
|  | 2 | Present | = 0 | 0.9438 | 0.9456 | 0.9460 | 0.9461 | 0.9461 | 0.9462 |
|  |  | SSDT [60] | ≠ 0 | 0.9233 | 0.9233 | 0.9233 | 0.9233 | 0.9233 | 0.9233 |
|  |  | HSDT[53] | ≠ 0 | 0.8406 | 0.8406 | 0.8406 | - | - | - |
|  | 5 | Present | = 0 | 1.0282 | 1.0303 | 1.0309 | 1.0310 | 1.0311 | 1.0311 |
|  |  | SSDT [60] | ≠ 0 | 1.0010 | 1.0010 | 1.0010 | 1.0010 | 1.0010 | 1.0010 |
|  |  | HSDT[53] | ≠ 0 | 0.9676 | 0.9676 | 0.9676 | - | - | - |
|  | 10 | Present | = 0 | 1.0801 | 1.0824 | 1.0830 | 1.0832 | 1.0832 | 1.0832 |
|  |  | SSDT [60] | ≠ 0 | 1.0487 | 1.0487 | 1.0487 | 1.0487 | 1.0487 | 1.0487 |
|  |  | HSDT[53] | ≠ 0 | 1.0879 | 1.0879 | 1.0879 | - | - | - |
| 10 | 0 | Present | = 0 | 0.7577 | 0.7608 | 0.7616 | 0.7618 | 0.7618 | 0.7618 |
|  |  | SSDT [60] | ≠ 0 | 0.7561 | 0.7561 | 0.7561 | 0.7561 | 0.7561 | 0.7561 |
|  | 1 | Present | = 0 | 0.8769 | 0.8821 | 0.8835 | 0.8838 | 0.8839 | 0.8839 |
|  |  | SSDT [60] | ≠ 0 | 0.8754 | 0.8754 | 0.8754 | 0.8754 | 0.8754 | 0.8754 |
|  | 2 | Present | = 0 | 0.9398 | 0.9464 | 0.9481 | 0.9486 | 0.9486 | 0.9487 |
|  |  | SSDT [60] | ≠ 0 | 0.9373 | 0.9373 | 0.9373 | 0.9373 | 0.9373 | 0.9373 |
|  | 5 | Present | = 0 | 1.0226 | 1.0310 | 1.0331 | 1.0336 | 1.0337 | 1.0338 |
|  |  | SSDT [60] | ≠ 0 | 1.0170 | 1.0170 | 1.0170 | 1.0170 | 1.0170 | 1.0170 |
|  | 10 | Present | = 0 | 1.0739 | 1.0829 | 1.0852 | 1.0859 | 1.0860 | 1.0860 |
|  |  | SSDT [60] | ≠ 0 | 1.0658 | 1.0658 | 1.0658 | 1.0658 | 1.0658 | 1.0658 |

**Table 14** Non-dimensional axial stress of (1-1-1) FG sandwich SS curved beams (Type B).

| L/h | p | Model | $\varepsilon_z$ | R/L | | | | | |
|---|---|---|---|---|---|---|---|---|---|
| | | | | 5 | 10 | 20 | 50 | 100 | ∞ |
| 5 | 0 | Present | = 0 | 3.8085 | 3.8241 | 3.8306 | 3.8340 | 3.8351 | 3.8362 |
| | | SSDT [60] | ≠ 0 | 3.8221 | 3.8402 | 3.8492 | 3.8546 | 3.8564 | 3.8582 |
| | | HSDT[53] | ≠ 0 | 3.7557 | 3.7557 | 3.7557 | - | - | - |
| | 1 | Present | = 0 | 1.4387 | 1.4454 | 1.4480 | 1.4494 | 1.4498 | 1.4502 |
| | | SSDT [60] | ≠ 0 | 1.4449 | 1.4526 | 1.4564 | 1.4587 | 1.4594 | 1.4602 |
| | | HSDT[53] | ≠ 0 | 1.5068 | 1.5068 | 1.5068 | - | - | - |
| | 2 | Present | = 0 | 1.9425 | 1.9521 | 1.9558 | 1.9577 | 1.9583 | 1.9588 |
| | | SSDT [60] | ≠ 0 | 1.9518 | 1.9625 | 1.9679 | 1.9711 | 1.9721 | 1.9732 |
| | | HSDT[53] | ≠ 0 | 1.9403 | 1.9403 | 1.9403 | - | - | - |
| | 5 | Present | = 0 | 2.6169 | 2.6307 | 2.6359 | 2.6384 | 2.6392 | 2.6399 |
| | | SSDT [60] | ≠ 0 | 2.6274 | 2.6424 | 2.6499 | 2.6544 | 2.6559 | 2.6575 |
| | | HSDT[53] | ≠ 0 | 2.5838 | 2.5838 | 2.5838 | - | - | - |
| | 10 | Present | = 0 | 2.9341 | 2.9498 | 2.9556 | 2.9585 | 2.9594 | 2.9602 |
| | | SSDT [60] | ≠ 0 | 2.9432 | 2.9604 | 2.9690 | 2.9741 | 2.9759 | 2.9776 |
| | | HSDT[53] | ≠ 0 | 2.8837 | 2.8837 | 2.8837 | - | - | - |
| 10 | 0 | Present | = 0 | 7.4633 | 7.5400 | 7.5689 | 7.5834 | 7.5877 | 7.5917 |
| | | SSDT [60] | ≠ 0 | 7.5665 | 7.6059 | 7.6256 | 7.6374 | 7.6414 | 7.6453 |
| | 1 | Present | = 0 | 2.8200 | 2.8558 | 2.8685 | 2.8745 | 2.8762 | 2.8778 |
| | | SSDT [60] | ≠ 0 | 2.8708 | 2.8875 | 2.8958 | 2.9008 | 2.9025 | 2.9041 |
| | 2 | Present | = 0 | 3.8072 | 3.8606 | 3.8792 | 3.8876 | 3.8899 | 3.8921 |
| | | SSDT [60] | ≠ 0 | 3.8840 | 3.9074 | 3.9191 | 3.9261 | 3.9284 | 3.9308 |
| | 5 | Present | = 0 | 5.1290 | 5.2079 | 5.2347 | 5.2465 | 5.2498 | 5.2527 |
| | | SSDT [60] | ≠ 0 | 5.2363 | 5.2691 | 5.2855 | 5.2954 | 5.2987 | 5.3020 |
| | 10 | Present | = 0 | 5.7516 | 5.8423 | 5.8728 | 5.8862 | 5.8899 | 5.8931 |
| | | SSDT [60] | ≠ 0 | 5.8687 | 5.9062 | 5.9250 | 5.9362 | 5.9400 | 5.9437 |

**Table 15.** Non-dimensional vertical displacement of (1-1-1) FG sandwich CC and CF curved beams (Type B).

| BC | L/h | p | R/L | | | | | |
|---|---|---|---|---|---|---|---|---|
| | | | 5 | 10 | 20 | 50 | 100 | ∞ |
| CC | 5 | 0 | 0.8170 | 0.8327 | 0.8368 | 0.8379 | 0.8381 | 0.8381 |
| | | 1 | 1.4579 | 1.4940 | 1.5033 | 1.5059 | 1.5063 | 1.5064 |
| | | 2 | 1.8816 | 1.9340 | 1.9476 | 1.9514 | 1.9520 | 1.9522 |
| | | 5 | 2.4408 | 2.5163 | 2.5359 | 2.5414 | 2.5422 | 2.5425 |
| | | 10 | 2.7063 | 2.7920 | 2.8143 | 2.8206 | 2.8215 | 2.8218 |
| | 10 | 0 | 0.5965 | 0.6300 | 0.6390 | 0.6415 | 0.6419 | 0.6420 |
| | | 1 | 1.1409 | 1.2314 | 1.2563 | 1.2635 | 1.2645 | 1.2648 |
| | | 2 | 1.4996 | 1.6377 | 1.6764 | 1.6875 | 1.6891 | 1.6896 |
| | | 5 | 1.9719 | 2.1789 | 2.2376 | 2.2546 | 2.2571 | 2.2579 |
| | | 10 | 2.1965 | 2.4346 | 2.5024 | 2.5220 | 2.5249 | 2.5258 |
| CF | 5 | 0 | 28.6872 | 28.7203 | 28.7286 | 28.7309 | 28.7312 | 28.7313 |
| | | 1 | 58.0282 | 58.1279 | 58.1529 | 58.1599 | 58.1609 | 58.1612 |
| | | 2 | 78.1225 | 78.2811 | 78.3209 | 78.3321 | 78.3337 | 78.3342 |
| | | 5 | 105.048 | 105.295 | 105.357 | 105.374 | 105.377 | 105.377 |
| | | 10 | 117.735 | 118.022 | 118.094 | 118.114 | 118.117 | 118.118 |
| | 10 | 0 | 27.7389 | 27.8655 | 27.8973 | 27.9062 | 27.9075 | 27.9079 |
| | | 1 | 56.6357 | 57.0216 | 57.1190 | 57.1463 | 57.1502 | 57.1515 |
| | | 2 | 76.4136 | 77.0299 | 77.1856 | 77.2294 | 77.2356 | 77.2377 |
| | | 5 | 102.904 | 103.865 | 104.108 | 104.176 | 104.186 | 104.189 |
| | | 10 | 115.387 | 116.504 | 116.787 | 116.867 | 116.878 | 116.882 |

**Table 16.** Non-dimensional vertical displacement of SS (2-2-1) FG sandwich curved beam (Type B).

| L/h | p | R/L | | | | | |
|---|---|---|---|---|---|---|---|
| | | 5 | 10 | 20 | 50 | 100 | ∞ |
| Transverse displacement ($\bar{w}$) | | | | | | | |
| 5 | 0 | 3.1652 | 3.1649 | 3.1651 | 3.1652 | 3.1652 | 3.1652 |
| | 1 | 5.8398 | 5.8380 | 5.8391 | 5.8396 | 5.8398 | 5.8399 |
| | 2 | 7.5567 | 7.5536 | 7.5556 | 7.5565 | 7.5567 | 7.5570 |
| | 5 | 9.7881 | 9.7829 | 9.7862 | 9.7877 | 9.7881 | 9.7885 |
| | 10 | 10.843 | 10.837 | 10.841 | 10.842 | 10.843 | 10.843 |
| 10 | 0 | 2.9453 | 2.9489 | 2.9498 | 2.9500 | 2.9501 | 2.9501 |
| | 1 | 5.5639 | 5.5755 | 5.5788 | 5.5801 | 5.5804 | 5.5806 |
| | 2 | 7.2502 | 7.2684 | 7.2740 | 7.2762 | 7.2768 | 7.2772 |
| | 5 | 9.4439 | 9.4717 | 9.4805 | 9.4842 | 9.4851 | 9.4860 |
| | 10 | 10.479 | 10.512 | 10.522 | 10.526 | 10.528 | 10.529 |
| Shear stress ($\bar{\tau}_{xz}$) | | | | | | | |
| 5 | 0 | 0.7534 | 0.7534 | 0.7534 | 0.7534 | 0.7534 | 0.7534 |
| | 1 | 0.8684 | 0.8682 | 0.8683 | 0.8684 | 0.8684 | 0.8684 |
| | 2 | 0.9274 | 0.9272 | 0.9274 | 0.9274 | 0.9274 | 0.9274 |
| | 5 | 1.0047 | 1.0044 | 1.0046 | 1.0047 | 1.0047 | 1.0047 |
| | 10 | 1.0516 | 1.0513 | 1.0515 | 1.0516 | 1.0516 | 1.0517 |
| 10 | 0 | 0.7595 | 0.7602 | 0.7604 | 0.7605 | 0.7605 | 0.7605 |
| | 1 | 0.8714 | 0.8728 | 0.8731 | 0.8733 | 0.8733 | 0.8733 |
| | 2 | 0.9295 | 0.9313 | 0.9318 | 0.9320 | 0.9320 | 0.9320 |
| | 5 | 1.0063 | 1.0084 | 1.0091 | 1.0093 | 1.0094 | 1.0094 |
| | 10 | 1.0531 | 1.0554 | 1.0561 | 1.0564 | 1.0565 | 1.0565 |
| Axial stress ($\bar{\sigma}_{xx}$) | | | | | | | |
| 5 | 0 | 3.8139 | 3.8158 | 3.8148 | 3.8141 | 3.8139 | 3.8136 |
| | 1 | 1.2517 | 1.2519 | 1.2518 | 1.2518 | 1.2517 | 1.2517 |
| | 2 | 1.5928 | 1.5927 | 1.5928 | 1.5928 | 1.5928 | 1.5928 |
| | 5 | 2.0260 | 2.0254 | 2.0258 | 2.0260 | 2.0260 | 2.0261 |
| | 10 | 2.2270 | 2.2262 | 2.2267 | 2.2270 | 2.2270 | 2.2271 |
| 10 | 0 | 7.5538 | 7.5533 | 7.5507 | 7.5483 | 7.5474 | 7.5465 |
| | 1 | 2.4786 | 2.4813 | 2.4816 | 2.4815 | 2.4814 | 2.4812 |
| | 2 | 3.1534 | 3.1590 | 3.1602 | 3.1604 | 3.1604 | 3.1604 |
| | 5 | 4.0105 | 4.0204 | 4.0232 | 4.0241 | 4.0243 | 4.0244 |
| | 10 | 4.4083 | 4.4202 | 4.4237 | 4.4250 | 4.4253 | 4.4255 |

**Table 17.** Non-dimensional vertical displacement of (1-8-1) FG sandwich straight beams (Type C) with different boundary condition.

| L/h | BC | Reference | Model | p= 0 | p = 1 | p= 2 | p = 5 | p = 10 |
|---|---|---|---|---|---|---|---|---|
| 5 | SS | Present | FE-PSDT | 3.9551 | 6.7126 | 8.0039 | 9.0717 | 9.4872 |
| | | Vo et al. [46] | Navier | 3.9788 | 6.7166 | 8.0083 | 9.0691 | 9.4817 |
| | | Vo et al. [46] | FE-TBT | 3.9788 | 6.7166 | 8.0083 | 9.0691 | 9.4817 |
| | CC | Present | FE-PSDT | 1.0093 | 1.6841 | 2.0524 | 2.5036 | 2.7466 |
| | | Vo et al. [46] | FE-TBT | 1.0273 | 1.7079 | 2.0825 | 2.5386 | 2.7866 |
| | CF | Present | FE-PSDT | 36.216 | 61.681 | 73.175 | 81.469 | 84.148 |
| | | Vo et al. [46] | FE-TBT | 36.468 | 61.737 | 73.244 | 81.533 | 84.216 |
| 20 | SS | Present | FE-PSDT | 3.6697 | 6.2602 | 7.4029 | 8.1531 | 8.3571 |
| | | Vo et al. [46] | Navier | 3.6934 | 6.2638 | 7.4085 | 8.1587 | 8.3434 |
| | | Vo et al. [46] | FE-TBT | 3.6934 | 6.2638 | 7.4085 | 8.1587 | 8.3619 |
| | CC | Present | FE-PSDT | 0.7476 | 1.2706 | 1.5044 | 1.6693 | 1.7210 |
| | | Vo et al. [46] | FE-TBT | 0.7536 | 1.2759 | 1.5122 | 1.6784 | 1.7300 |
| | CF | Present | FE-PSDT | 35.121 | 59.948 | 70.875 | 77.958 | 79.831 |
| | | Vo et al. [46] | FE-TBT | 35.349 | 59.966 | 70.901 | 77.988 | 79.858 |

**Table 18.** Non-dimensional axial stress and transverse shear stress of (1-8-1) FG sandwich SS straight beams (Type C).

| p | Reference | Model | L/h = 5 | | L/h = 20 | |
|---|---|---|---|---|---|---|
| | | | $\bar{\sigma}_{xx}$ | $\bar{\tau}_{xz}$ | $\bar{\sigma}_{xx}$ | $\bar{\tau}_{xz}$ |
| 0 | Present | FE-PSDT | 4.4610 | 0.7802 | 17.6144 | 0.7878 |
| | Vo et al. [46] | Navier | 4.4636 | 0.7597 | 17.6327 | 0.7702 |
| | Vo et al. [46] | FE-TBT | 4.4660 | 0.7611 | 17.6400 | 0.7785 |
| 1 | Present | FE-PSDT | 6.0312 | 0.7519 | 23.7834 | 0.7610 |
| | Vo et al. [46] | Navier | 6.0094 | 0.7318 | 23.7080 | 0.7436 |
| | Vo et al. [46] | FE-TBT | 6.0130 | 0.7315 | 23.7200 | 0.7416 |
| 2 | Present | FE-PSDT | 6.5497 | 0.6647 | 25.7654 | 0.6738 |
| | Vo et al. [46] | Navier | 6.5256 | 0.6445 | 25.6849 | 0.6558 |
| | Vo et al. [46] | FE-TBT | 6.5290 | 0.6432 | 25.7000 | 0.6452 |
| 5 | Present | FE-PSDT | 6.9186 | 0.5527 | 27.0641 | 0.5619 |
| | Vo et al. [46] | Navier | 6.8886 | 0.5319 | 26.9694 | 0.5425 |
| | Vo et al. [46] | FE-TBT | 6.8930 | 0.5316 | 269800 | 0.5400 |
| 10 | Present | FE-PSDT | 7.2565 | 0.6021 | 28.3354 | 0.6125 |
| | Vo et al. [46] | Navier | 7.2229 | 0.5792 | 28.2283 | 0.5910 |
| | Vo et al. [46] | FE-TBT | 7.2270 | 0.5798 | 28.2400 | 0.5969 |

**Table 19.** Non-dimensional transverse displacement and stresses of (1-8-1) FG sandwich SS curved beams (Type C).

| L/h | p | R/L | | | | | |
|---|---|---|---|---|---|---|---|
| | | 5 | 10 | 20 | 50 | 100 | ∞ |
| Transverse displacement ($\bar{w}$) | | | | | | | |
| 5 | 0 | 3.9518 | 3.9540 | 3.9547 | 3.9550 | 3.9551 | 3.9551 |
| | 1 | 6.7027 | 6.7083 | 6.7107 | 6.7119 | 6.7123 | 6.7126 |
| | 2 | 7.9918 | 7.9986 | 7.9986 | 8.0030 | 8.0034 | 8.0039 |
| | 5 | 9.0591 | 9.0662 | 9.0691 | 9.0707 | 9.0712 | 9.0717 |
| | 10 | 9.4752 | 9.4820 | 9.4848 | 9.4863 | 9.4868 | 9.4872 |
| 10 | 0 | 3.7173 | 3.7238 | 3.7257 | 3.7265 | 3.7267 | 3.7268 |
| | 1 | 6.3271 | 6.3415 | 6.3468 | 6.3493 | 6.3501 | 6.3508 |
| | 2 | 7.4955 | 7.5120 | 7.5183 | 7.5214 | 7.5224 | 7.5232 |
| | 5 | 8.3094 | 8.3257 | 8.3320 | 8.3352 | 8.3361 | 8.3370 |
| | 10 | 8.5573 | 8.5728 | 8.5787 | 8.5817 | 8.5826 | 8.5834 |
| Shear stress ($\bar{\tau}_{xz}$) | | | | | | | |
| 5 | 0 | 0.7798 | 0.7801 | 0.7802 | 0.7802 | 0.7802 | 0.7802 |
| | 1 | 0.7513 | 0.7517 | 0.7518 | 0.7519 | 0.7519 | 0.7519 |
| | 2 | 0.6641 | 0.6644 | 0.6644 | 0.6646 | 0.6647 | 0.6647 |
| | 5 | 0.5523 | 0.5525 | 0.5526 | 0.5527 | 0.5527 | 0.5527 |
| | 10 | 0.6016 | 0.6019 | 0.6020 | 0.6021 | 0.6021 | 0.6021 |
| 10 | 0 | 0.7847 | 0.7857 | 0.7860 | 0.7861 | 0.7861 | 0.7861 |
| | 1 | 0.7572 | 0.7583 | 0.7587 | 0.7589 | 0.7589 | 0.7590 |
| | 2 | 0.6702 | 0.6711 | 0.6715 | 0.6716 | 0.6717 | 0.6717 |
| | 5 | 0.5586 | 0.5593 | 0.5596 | 0.5597 | 0.5597 | 0.5597 |
| | 10 | 0.6089 | 0.6096 | 0.6098 | 0.6099 | 0.6100 | 0.6100 |
| Axial stress ($\bar{\sigma}_{xx}$) | | | | | | | |
| 5 | 0 | 4.4619 | 4.4620 | 4.4617 | 4.4613 | 4.4612 | 4.4610 |
| | 1 | 6.0229 | 6.0277 | 6.0296 | 6.0306 | 6.0309 | 6.0312 |
| | 2 | 6.5393 | 6.5451 | 6.5451 | 6.5489 | 6.5493 | 6.5497 |
| | 5 | 6.9094 | 6.9146 | 6.9167 | 6.9179 | 6.9182 | 6.9186 |
| | 10 | 7.2492 | 7.2534 | 7.2551 | 7.2560 | 7.2562 | 7.2565 |
| 10 | 0 | 8.8254 | 8.8320 | 8.8321 | 8.8312 | 8.8307 | 8.8301 |
| | 1 | 11.8842 | 11.9099 | 11.9191 | 11.9235 | 11.9247 | 11.9258 |
| | 2 | 12.8768 | 12.9061 | 12.9172 | 12.9228 | 12.9245 | 12.9260 |
| | 5 | 13.5499 | 13.5757 | 13.5854 | 13.5903 | 13.5917 | 13.5931 |
| | 10 | 14.2012 | 14.2232 | 14.2310 | 14.2347 | 14.2358 | 14.2368 |

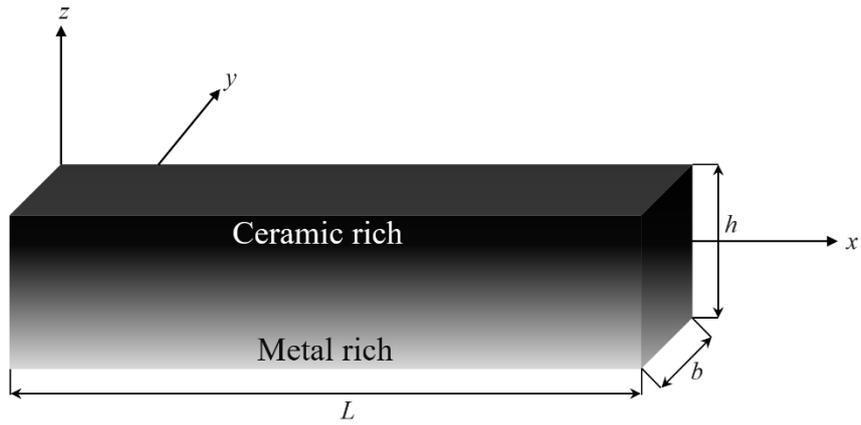

(a) Coordinate system of FG beam

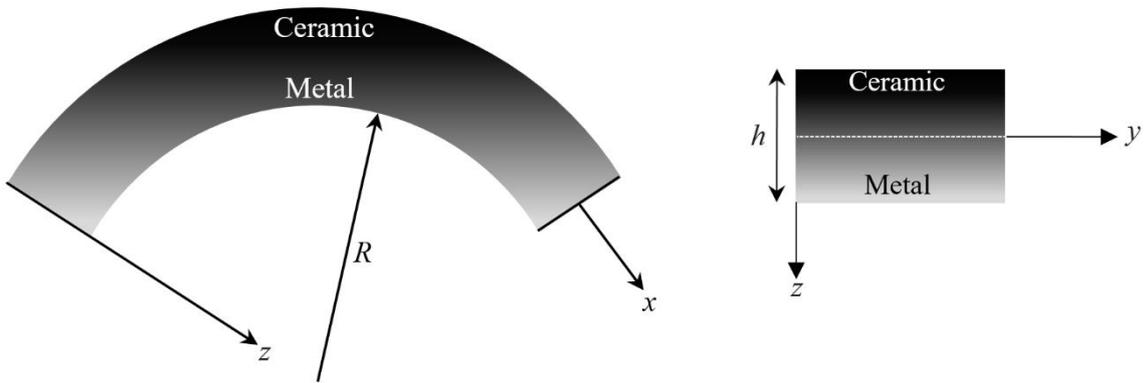

(b) FG beam (Type A)

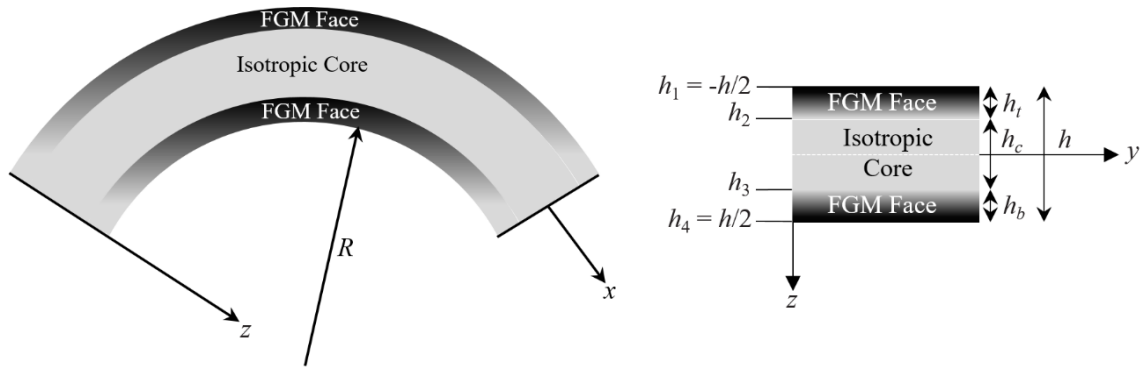

(c) Sandwich beam with FG face sheets and isotropic core (Type B)

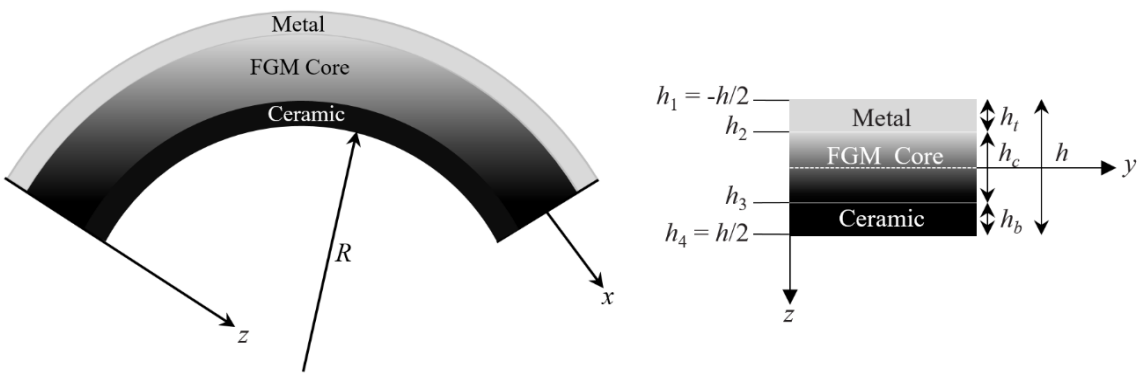

(d) Sandwich beam with isotropic face sheets and FG core (Type C)

**Fig. 1** Geometry and coordinate of isotropic and FG sandwich curved beams.

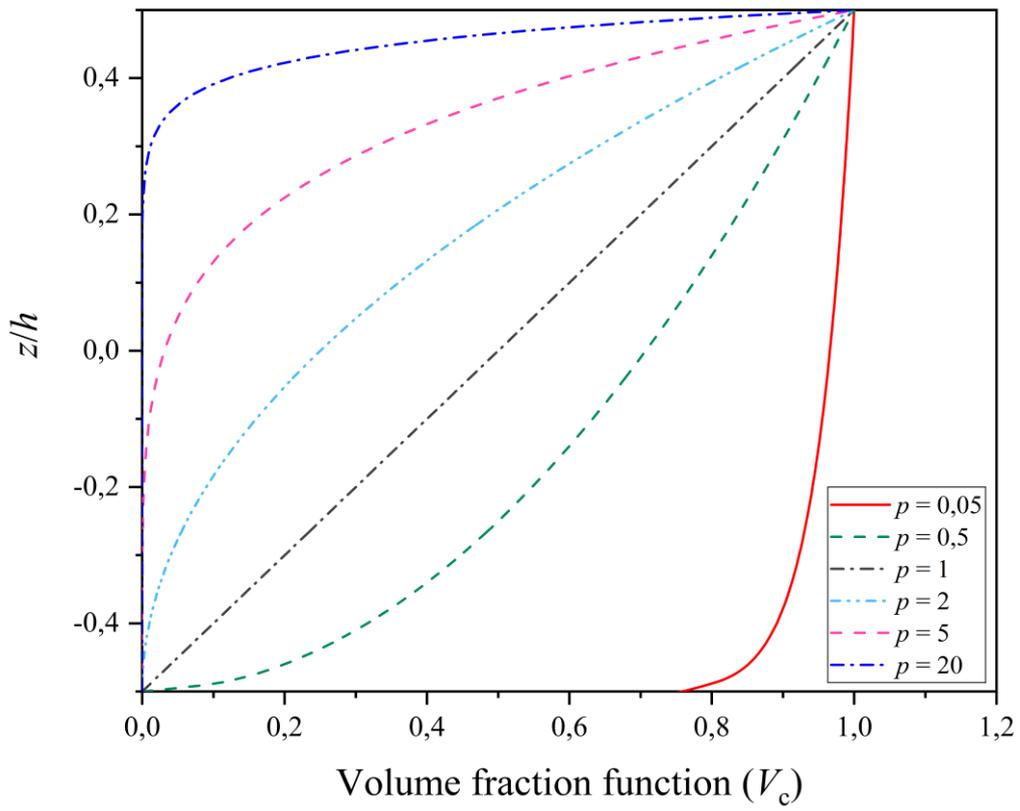

(a) Type (A)

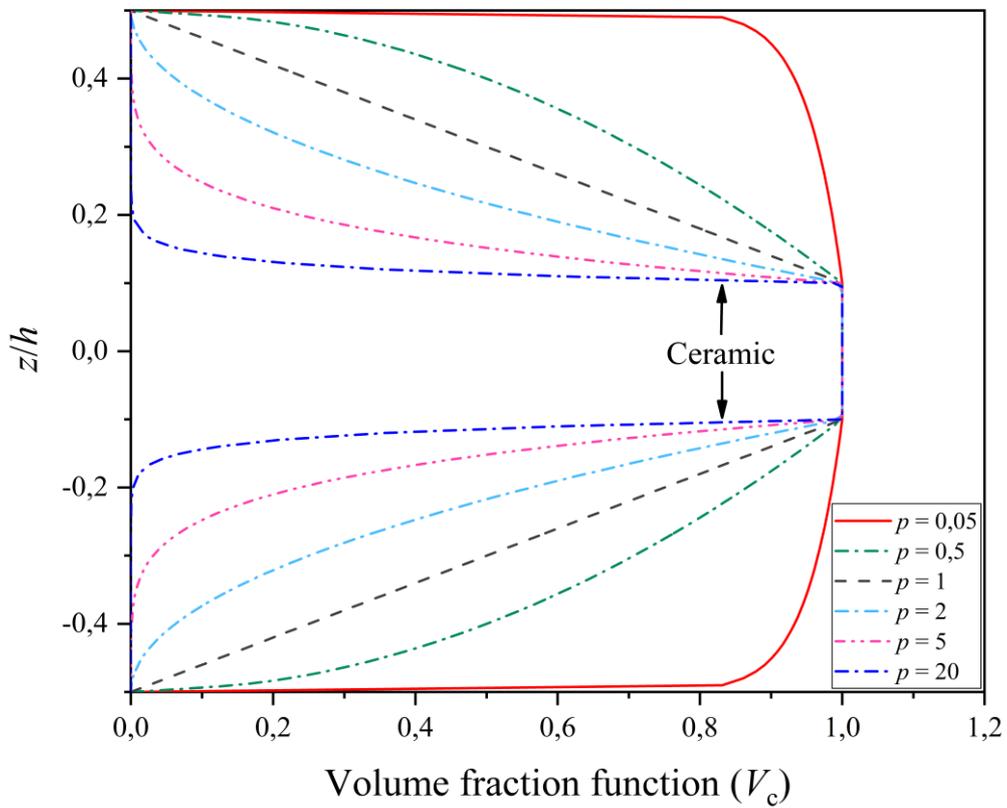

(b) Type (B) 2-1-2

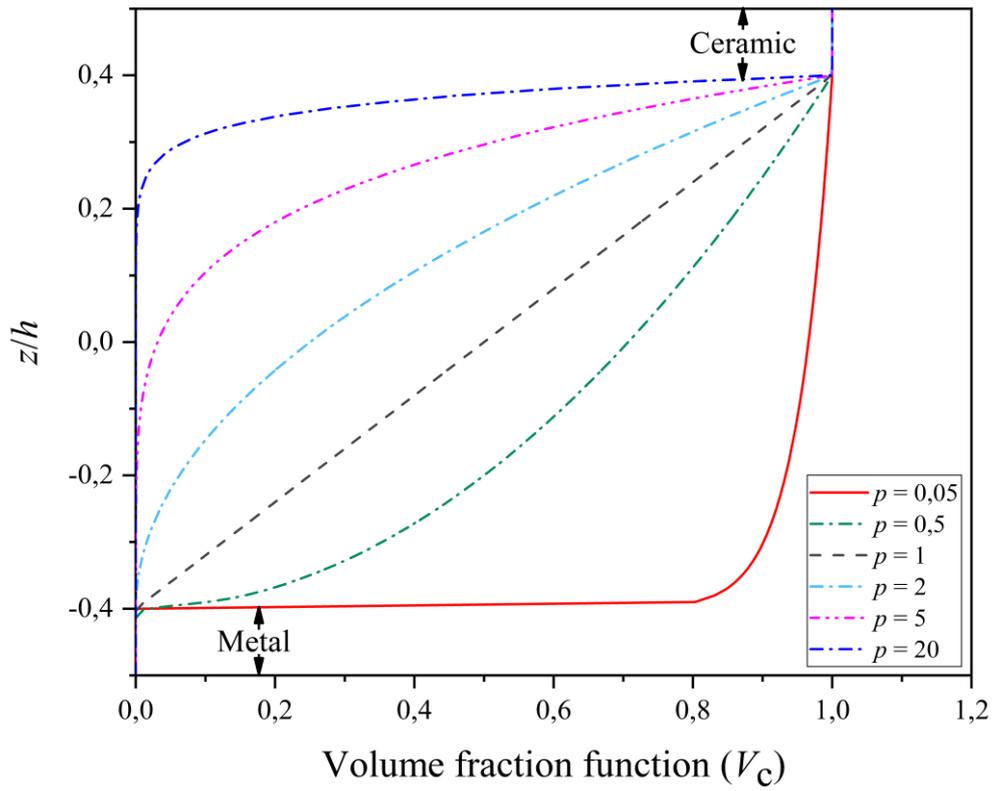

(c) Type (C) 1-8-1

**Fig. 2** Variation of the volume fraction function through the thickness of three types of FGM beams for various values of the power law index $p$.

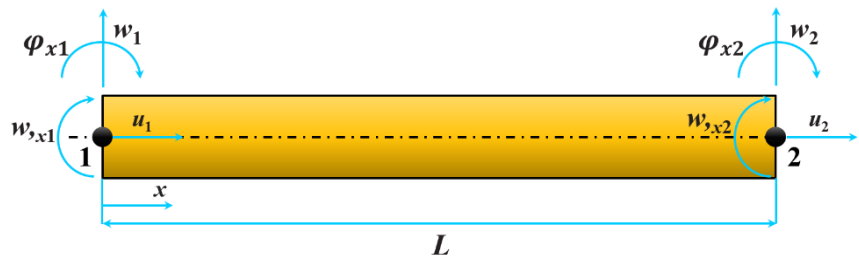

**Fig. 3** Present two-noded beam element with corresponding DOFs.

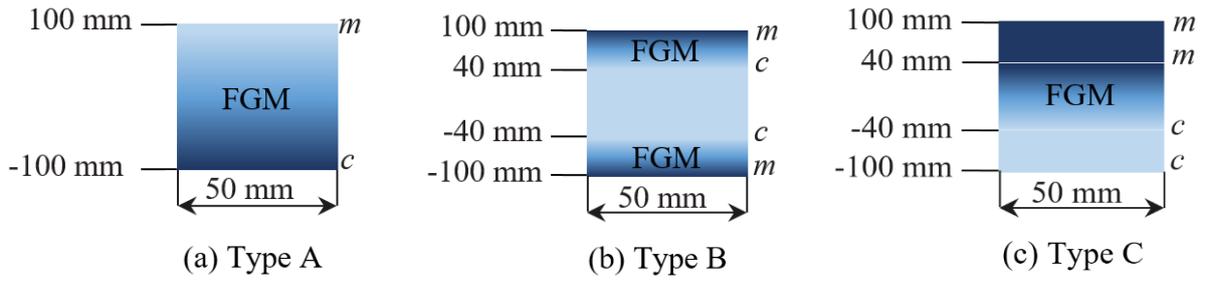

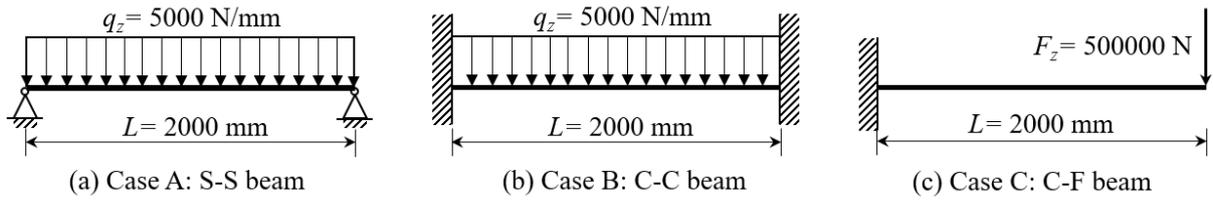

**Fig.4** Studied FG beams with different boundary conditions, material distribution and loads.

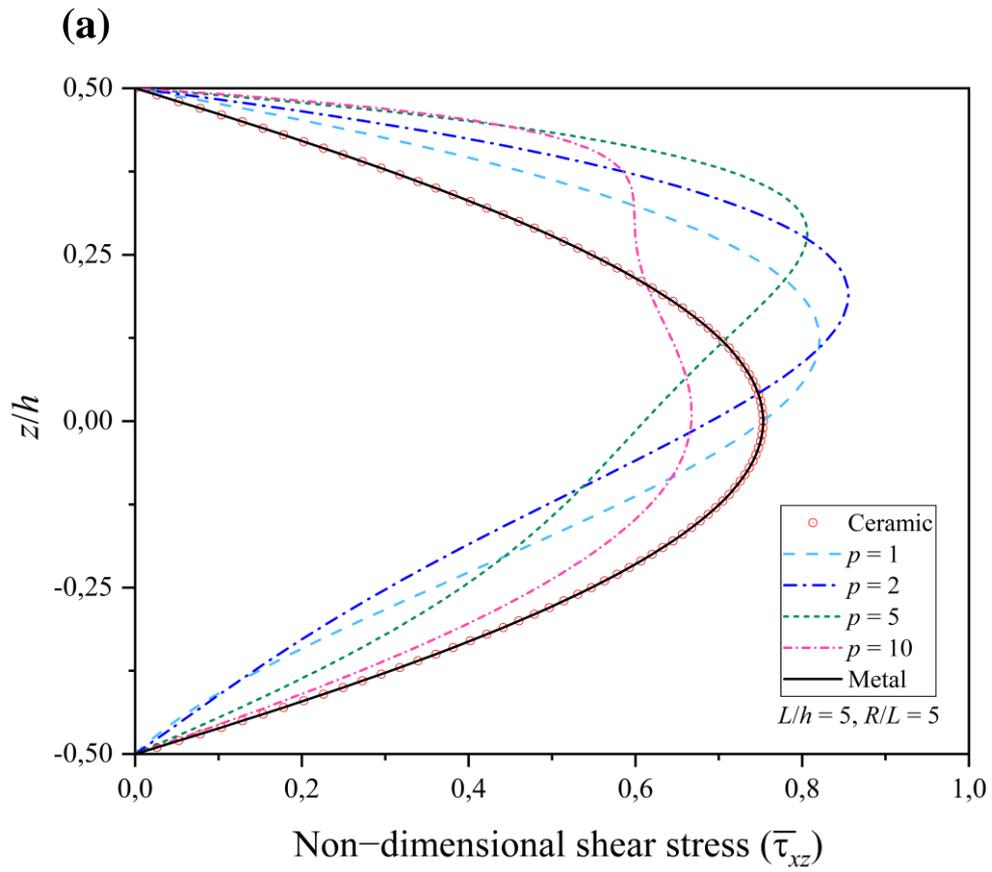

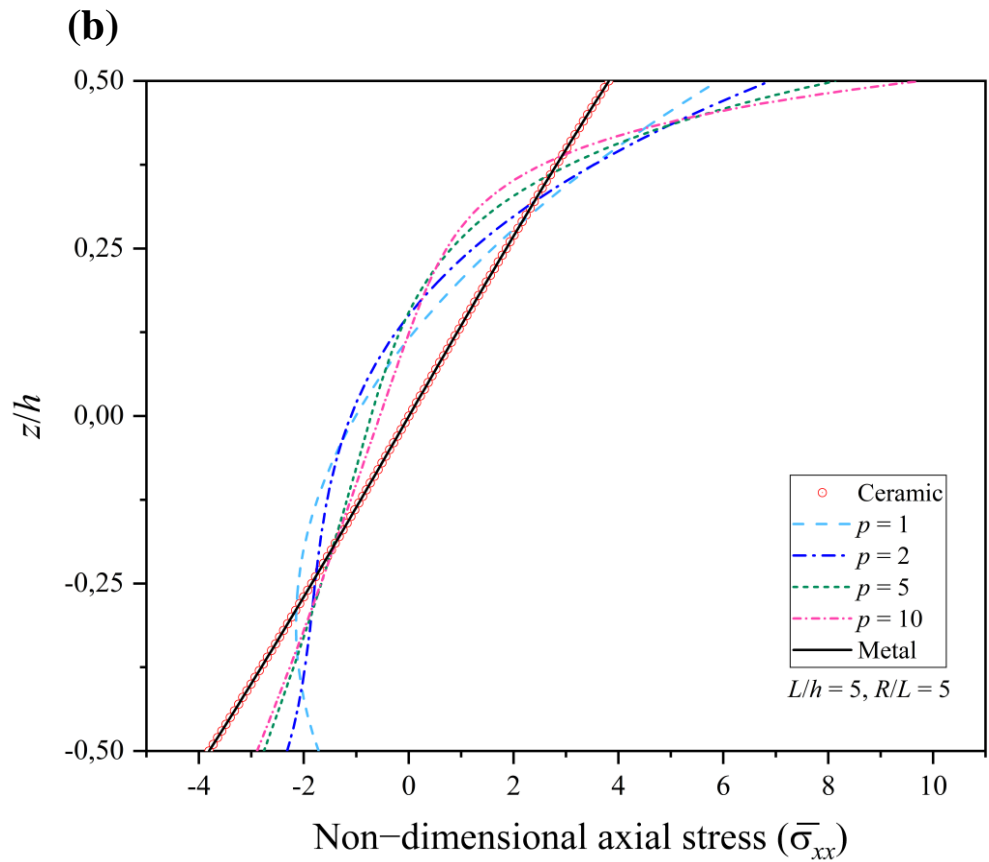

**Fig. 5** Distribution of non-dimensional stresses through the thickness of single layer FG SS curved beam (Type A), (a) transverse shear stress, (b) axial stress.

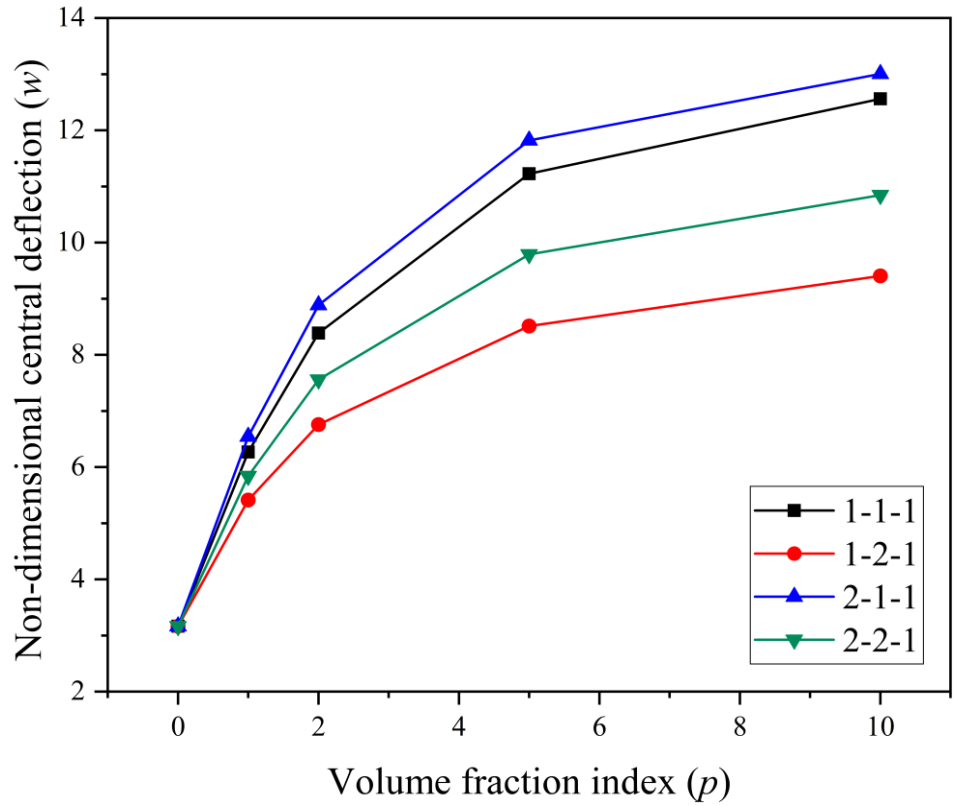

**Fig. 6** Effect of volume fraction index (*p*) with different core-to-face sheets thickness ratio on the non-dimensional center deflection of SS sandwich straight beams with FG face sheets ($a/h = 5$).

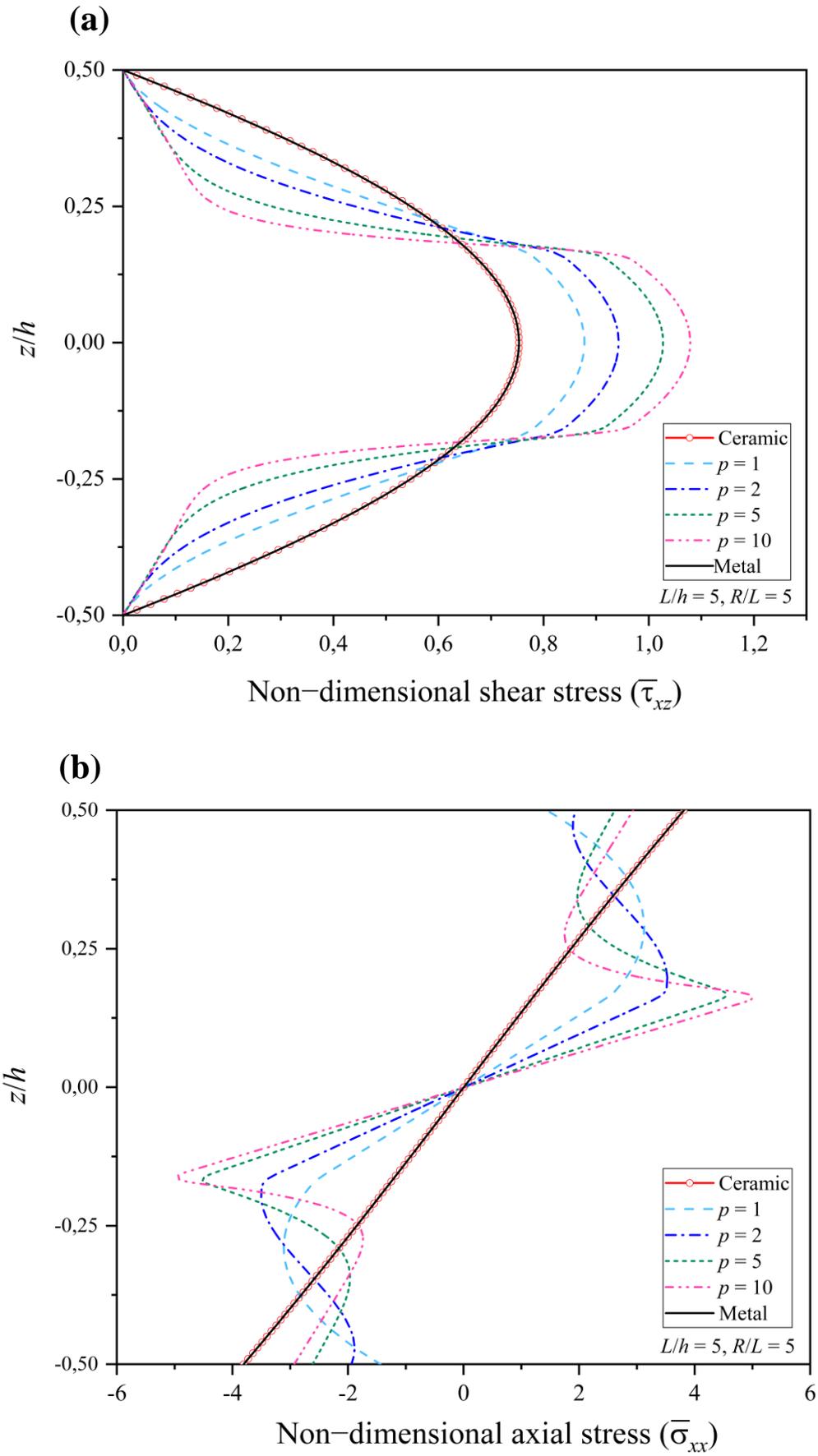

**Fig. 7** Distribution of non-dimensional stresses along the thickness of symmetric (1-1-1) FG sandwich SS curved beams (Type B), (a) transverse shear stress, (b) axial stress.

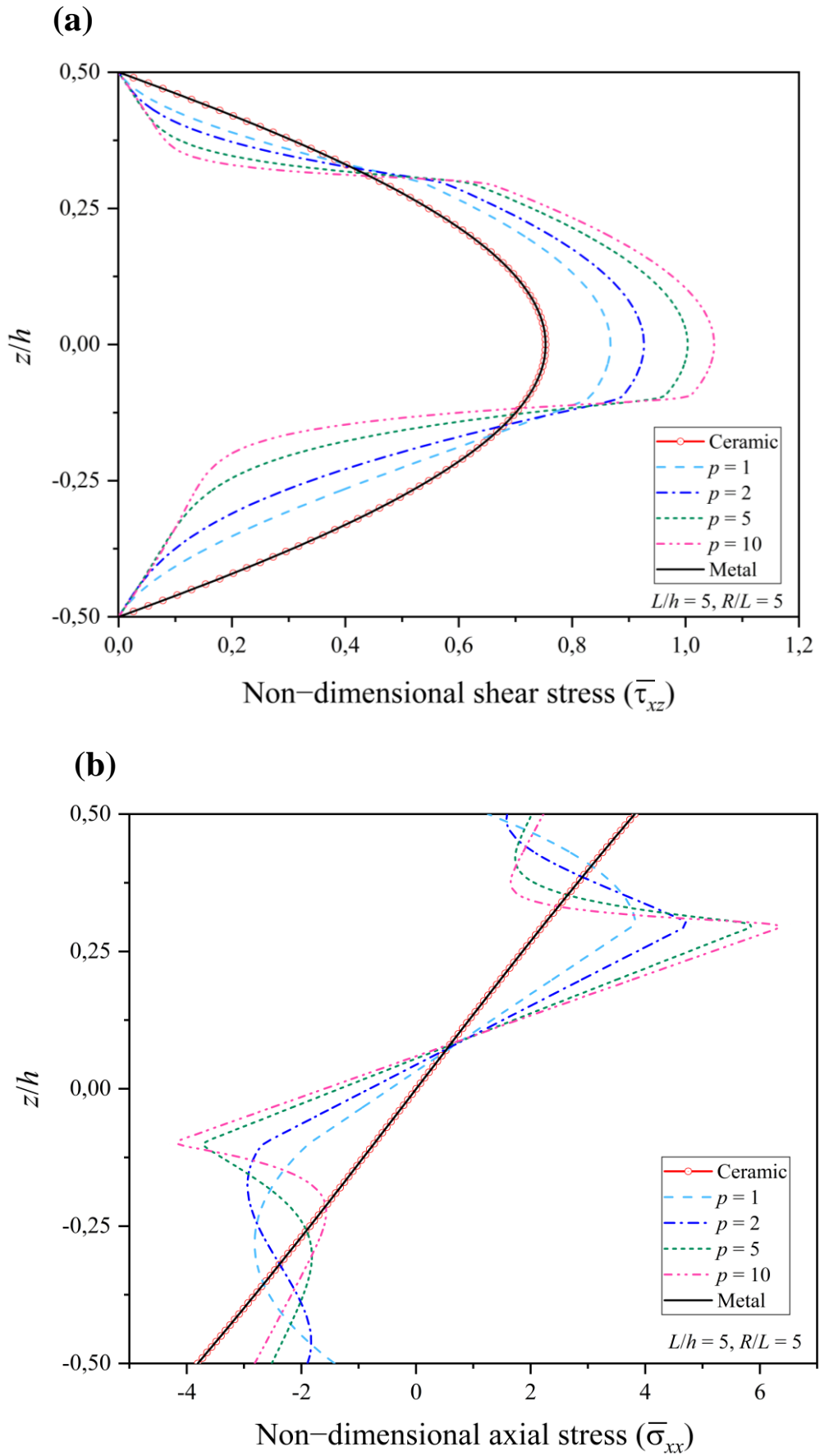

**Fig. 8** Distribution of non-dimensional stresses along the thickness of non-symmetric (2-2-1) FG sandwich SS curved beams (Type B), (a) transverse shear stress, (b) axial stress.

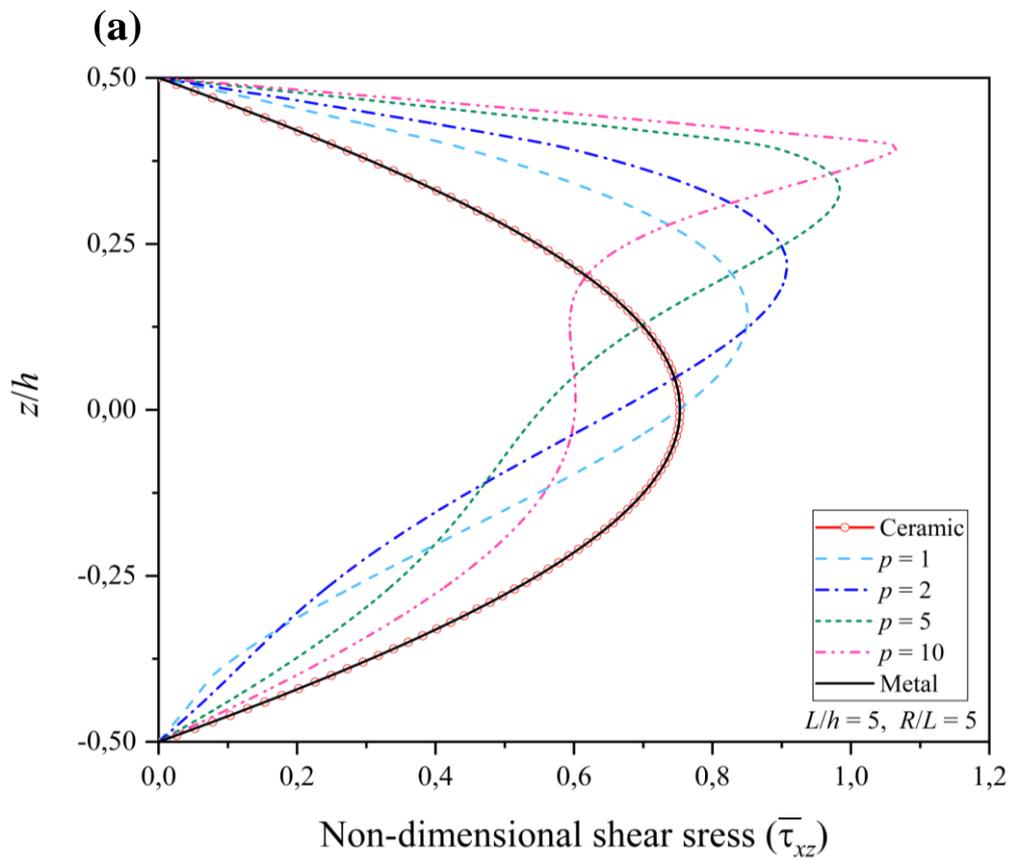

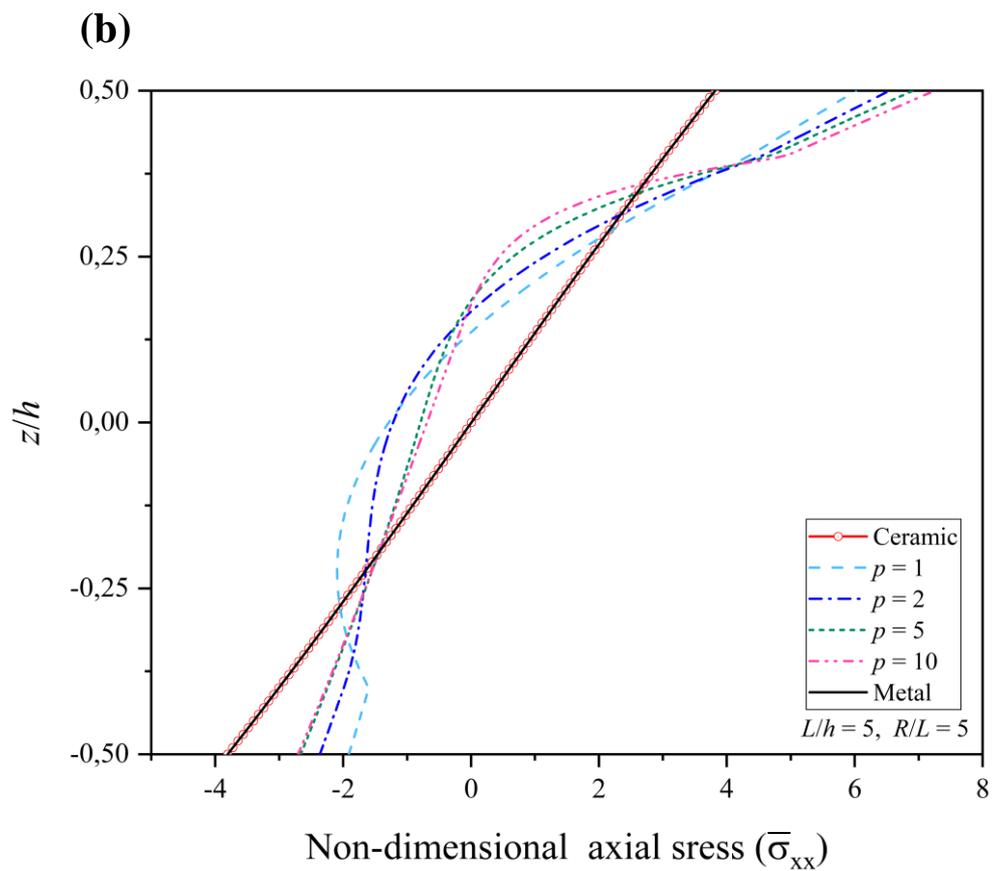

**Fig. 9** Distribution of non-dimensional stresses along the thickness of (1-8-1) FG sandwich SS curved beams (Type C), (a) transverse shear stress, (b) axial stress.